\PassOptionsToPackage{tracking=false}{microtype}
\documentclass[11pt, a4paper, logo]{googledeepmind}
\usepackage[dvipsnames]{xcolor}         % colors
\usepackage{amsmath}
\usepackage{mathtools}
\usepackage{cprotect}
\mathtoolsset{showonlyrefs}
\usepackage{amssymb}
\usepackage{bm}
\usepackage{subcaption}
\usepackage{wrapfig}
\usepackage{enumitem}
\usepackage{tikz}
\usetikzlibrary{arrows.meta, calc, decorations.pathmorphing, shapes, backgrounds}
\theoremstyle{plain}

\newtheorem{theorem}{Theorem}[section]
\newtheorem{proposition}[theorem]{Proposition}

\theoremstyle{definition}
\newtheorem{definition}[theorem]{Definition}

\theoremstyle{remark}

\usepackage{algorithm}
\usepackage[noend]{algpseudocode}

\newcommand{\Z}{\mathbb{Z}}
\newcommand{\E}{\mathbf{E}}

\newcommand{\Sn}{\mathbb{S}_n}

\newcommand{\cS}{\mathcal{S}}
\newcommand{\cL}{\mathcal{L}}
\newcommand{\cQ}{\mathcal{Q}}

\newcommand{\fG}{\mathfrak{g}}

\newcommand{\ones}{\mathbf{1}}

\newcommand{\ud}{\mathrm{d}}

\newcommand{\GL}{\text{GL}}
\newcommand{\SL}{\text{SL}}
\newcommand{\SO}{\text{SO}}

\newcommand{\Or}{\text{O}}
\newcommand{\Ad}{\text{Ad}}
\newcommand{\gA}{\mathfrak{g}}

\definecolor{CFBlue}{cmyk}{0.7,0.3,0,0.1}
\colorlet{yh}{CFBlue}


\usepackage{enumitem}   
\usepackage{wrapfig}
\usepackage{float}
\usepackage{url}
\usepackage{upgreek}
\usepackage{tikz-cd}

\newcommand{\R}{\mathbb R}
\newcommand{\RR}{\mathbb R}

\newcommand{\cN}{\mathcal N}
\newcommand{\cM}{\mathcal M}

\newcommand{\cX}{\mathcal X}

\newcommand{\cP}{\mathcal{P}}
\newcommand{\KL}{\mathop{\mathrm{KL}}\nolimits}

\renewcommand{\top}{\intercal}

\usepackage{nicefrac}

\usepackage[authoryear, sort&compress, round]{natbib}
\title{Statistical Geometry of Group Actions: Symmetry, Transport, Estimation}
\author{Q. Berthet, P. Bartlett, R. Elie, Y.-H. Wu}
\date{\today}
\begin{abstract}
We study variational inference in the presence of distributional symmetry, focusing on the interaction between two distinct group actions on a sample space: a \emph{symmetry} action that leaves a reference distribution invariant, and a \emph{transport} action that generates a parametric variational family by pushforward. By analyzing the algebraic structure of the total transformation group generated by both actions and the induced adjoint representation of the symmetry group on the transport Lie algebra, we establish that the Riemannian gradient of the Kullback--Leibler (KL) divergence at the true transport parameter is constrained to lie in the trivial representation of the symmetry group. Consequently, symmetry forces the gradient to vanish along all non-trivial representations. Combining this representation-theoretic constraint with parameter identifiability yields an orthogonal decomposition of the transport Lie algebra into three interpretable subspaces: an \emph{identifiability kernel} (along which the objective is locally constant), a \emph{bias subspace} (where the gradient may be nonzero due to model misspecification), and a \emph{symmetry-breaking subspace} (where the gradient is forced to vanish by symmetry alone). Finally, we characterize the residual bias directions in terms of distributional mismatch between the base model and the standardized target, and derive an explicit Stein-type moment-matching condition that eliminates this bias, extending the location-scale results of \citet{margossian2024variational} to general Lie group transformation families.
\end{abstract}
\microtypesetup{tracking=false}
\begin{document}
\maketitle

\section{Introduction}

In modern statistical practice, variational inference \citep{jordan1999introduction, wainwright2008graphical, blei2017variational} has become a primary methodology for scalable approximate inference, most notably for approximating intractable posterior distributions in Bayesian models. Rather than sampling directly from a complex target distribution $p$ on a sample space $\cX$, variational inference casts inference as an optimization problem: one seeks a distribution $q \in \cQ$ within a tractable \emph{variational family} that minimizes a statistical divergence---most commonly the Kullback--Leibler (KL) divergence $\KL(q \| p)$. In practice, the variational family $\cQ$ is nearly always misspecified with respect to $p$: factorized mean-field families sacrifice dependence structure for tractability, while structured transformation families and normalizing flows \citep{rezende2015variational, papamakarios2021normalizing} trade computational cost for richer geometry. Because $p$ rarely belongs to $\cQ$, a foundational statistical question arises: \emph{even under model misspecification, when does the variational optimum exactly recover meaningful structural parameters of the target distribution?}

Symmetry offers a natural geometric principle to address this question. In classical mathematical statistics, group invariance is foundational: equivariance principles govern the theory of optimal invariant estimators and hypothesis tests \citep{eaton1989group, wijsman1990invariant, lehmann2005testing}, the Hunt--Stein theorem establishes conditions under which minimax procedures can be chosen to be group-invariant \citep{bondar1981amenability}, and the interplay between sufficiency and invariance is central to statistical decision theory \citep{hall1965relationship}. In contemporary machine learning, group invariance and equivariance have similarly guided the design of structured neural architectures \citep{cohen2016group, weiler2019general, bronstein2021geometric}. By contrast, the role of distributional symmetry in shaping the optimization landscape of variational inference---specifically, how shared invariances between the variational base distribution and the standardized target distribution constrain stationary points of the variational objective---has only recently begun to receive systematic study.

A striking step in this direction was taken by \citet{margossian2024variational}, who demonstrated that variational inference can achieve exact parameter recovery in misspecified location families (generated by translations $x \mapsto x + \mu$) and location-scale families (generated by affine transformations $x \mapsto \mu + Ax$) whenever the base distribution $q_0$ and the standardized target share a common symmetry. Specifically, they established that the gradient of the KL divergence vanishes at the true location parameter under central reflection symmetry ($x \mapsto -x$), and at the true correlation matrix under elliptical symmetry, while uncovering a persistent scale bias when estimating the full covariance matrix. These findings led the authors to pose the following general question:
\begin{quote}
``\emph{Finally, an intriguing conjecture, also deserving of study, is that for each symmetry of $p$, there is a corresponding statistic that VI under mild conditions can exactly recover.}''
\end{quote}
In subsequent work, \citet{margossian2025generalized} extended these location-scale recovery guarantees from the standard variational objective $\KL(q \| p)$ to general $f$-divergences without requiring log-concavity. In the present paper, we address this conjecture in full generality by developing a Lie-theoretic and representation-theoretic framework for transformation-based variational inference. While our first-order symmetry constraint (Theorem~\ref{thm:gradient-constraint}) relies only on diffeomorphism invariance of the divergence and therefore applies broadly to any invariant loss functional, working explicitly with the KL divergence $\KL(q_g \| p)$ allows us to derive exact Stein-type moment-matching conditions that govern the remaining non-stationary directions. Our analysis reveals precisely how the algebraic compatibility between the symmetry group of the distribution and the Lie group of variational transports determines which parameter components are recovered without bias and which remain susceptible to model misspecification.

Concretely, we study the interaction between two group actions on the sample space $\cX$: a compact \emph{symmetry group} $H$ acting by isometries, and a \emph{transport Lie group} $G$ acting by diffeomorphisms to generate the variational family. This formulation connects classical statistical transformation models \citep{eaton1989group} and Lie group representation theory \citep{knapp2002lie, fulton2013representation, hall2015lie} with information geometry \citep{murray1993differential, amari2016information} and Riemannian optimization on manifolds \citep{absil2008optimization, boumal2023introduction}.

To formalize the estimation problem, we consider variational families obtained by pushing forward a fixed reference distribution $q_0 \in \cP(\cX)$ under the action of $G$: if $Z \sim q_0$, then $X = g \cdot Z \sim q_g$, yielding the variational optimization problem
\[
\min_{g \in G} \KL(q_g \| p)\, .
\]
On the data side, we assume that the target distribution $p$ is generated by transporting an unknown standardized distribution $\tilde p$ by a true parameter $g^* \in G$ that we wish to estimate; equivalently, $\tilde p$ is the \emph{pullback} of $p$ under $g^*$, meaning that if $X \sim p$, then the standardized variable $\xi = (g^*)^{-1} \cdot X$ follows $\tilde p$. Importantly, we do not assume that $\tilde p$ coincides with $q_0$, so the variational family $\cQ_G = \{q_g : g \in G\}$ may be misspecified ($p \notin \cQ_G$). Instead, we assume a shared \emph{structural symmetry}: both the variational base $q_0$ and the standardized target $\tilde p$ are invariant under the action of the symmetry group $H$. Geometrically, this means that $q_0$ and $\tilde p$ share the same notion of centering or alignment around the origin, and the statistical question is whether this shared invariance suffices to make $g^*$ a stationary point (and local minimizer) of the variational loss $g \mapsto \KL(q_g \| p)$.

\textbf{Main contributions.} We characterize the first-order geometry of transport variational inference at $g^*$ through the conjugation action of $H$ on $G$ and its linearization, the adjoint representation of $H$ on the Lie algebra $\gA = T_e G$:
\begin{enumerate}[leftmargin=*]
    \item \textbf{A Lie-theoretic framework for symmetric transport VI} (Sections~\ref{sec:problem}--\ref{sec:identifiability}). We formulate a general statistical framework pairing a compact symmetry group $H$ with an $H$-normalized transport Lie group $G$. Analyzing the total transformation group $K = GH$ and the intersection subgroup $N = G \cap H$ allows us to separate parameter identifiability on the coset manifold $\cM = G/N$ from first-order estimation bias, unifying location, scale, directional, and permutation models under a single geometric theory.
    
    \item \textbf{Gradient constraint and orthogonal Lie algebra decomposition} (Theorem~\ref{thm:gradient-constraint} and Proposition~\ref{prop:decomposition}). We prove that the Riemannian gradient $\nabla_G \cL(g^*)$ is constrained to lie in the trivial representation $\gA^H = \{u \in \gA : \Ad_h u = u \;\forall h \in H\}$ of the adjoint action of $H$ on $\gA$. Combining this invariance with the identifiability structure yields an orthogonal decomposition of the transport Lie algebra,
    \[
    \gA = \mathfrak{n} \;\oplus\; \gA_{\text{bias}}^H \;\oplus\; \gA^\perp\, ,
    \]
    into an \emph{identifiability kernel} $\mathfrak{n}$ (along which the loss is locally constant), a \emph{bias subspace} $\gA_{\text{bias}}^H = \gA^H \cap \mathfrak{n}^\perp$ (consisting of symmetry-preserving perturbations where model mismatch can induce a nonzero gradient), and a \emph{symmetry-breaking subspace} $\gA^\perp = (\gA^H)^\perp \cap \mathfrak{n}^\perp$ (along which the gradient is forced to vanish by symmetry alone). Whenever $\gA_{\text{bias}}^H = \{0\}$, $g^*$ is guaranteed to be a first-order stationary point regardless of model misspecification between $q_0$ and $\tilde p$.
   
    \item \textbf{Stein-type moment-matching condition and debiasing recipe} (Theorem~\ref{thm:moment-matching}). When the bias subspace $\gA_{\text{bias}}^H$ is nontrivial, we derive an exact necessary and sufficient condition for $g^*$ to be a critical point. Using integration by parts along the infinitesimal generators $V_u$ of the group action, we show that the gradient along $u \in \gA_{\text{bias}}^H$ vanishes if and only if the directional scores of the log-potentials of $q_0$ and $\tilde p$ match in expectation under $q_0$. This characterization provides a constructive calibration procedure to eliminate estimation bias by rescaling or shifting the base distribution $q_0$ along $\exp(\gA_{\text{bias}}^H)$ while preserving its $H$-invariance.
\end{enumerate}

\textbf{Organization.} The remainder of the paper is organized as follows. Section~\ref{sec:related} situates our work within the literature on variational inference, group invariance in statistics, and Riemannian optimization. Section~\ref{sec:problem} introduces transport variational families, $H$-conjugate symmetry, and the structural normalization condition, together with six guiding examples that we trace throughout the paper. Section~\ref{sec:identifiability} analyzes the total transformation group $K = GH$, the conjugation homomorphism, and the three resulting identifiability regimes on the quotient space $\cM = G/N$. Section~\ref{sec:lie} develops the required differential calculus on $G$ and introduces the adjoint representation $\Ad_H$ and its trivial subspace $\gA^H$. Section~\ref{sec:main-result} presents our main theoretical results: the gradient constraint theorem, the orthogonal decomposition of $\gA$, and the moment-matching theorem. Section~\ref{sec:discussion} discusses gauge-fixing vs.\ base calibration, second-order curvature conditions, and extensions. Technical proofs and detailed example derivations are deferred to the Appendices.

%% ===========================================================================
\section{Related work} \label{sec:related}
%% ===========================================================================

\textbf{Variational inference and transport-based families.} The variational formulation of approximate inference, rooted in statistical physics and convex duality, was brought to prominence in statistics and machine learning by \citet{jordan1999introduction}, \citet{wainwright2008graphical}, and \citet{blei2017variational}. The development of black-box variational inference and reparameterization gradients \citep{ranganath2014black, rezende2015variational} subsequently enabled optimization over general continuous families where expectations under $q_g$ are evaluated by pushing forward samples $Z \sim q_0$ through a differentiable transformation $x = g(Z)$. In the normalizing flow literature \citep{rezende2015variational, kobyzev2020normalizing, papamakarios2021normalizing}, the transformation $g$ is typically parameterized by a deep composition of general diffeomorphisms. Our framework studies the structured setting in which the transport maps form a finite-dimensional Lie group $G$ acting on $\cX$: this encompasses classical parametric families (such as location, scale, and directional families) while endowing the variational parameter space with an explicit Lie algebra that makes global and local symmetry analysis tractable.

\textbf{Group invariance in mathematical statistics and machine learning.} Group actions have long served as a unifying language in theoretical statistics \citep{eaton1989group, wijsman1990invariant, lehmann2005testing}. In classical statistical decision theory, when a family of distributions is closed under a group action, one seeks estimators that are \emph{equivariant}---such as the Pitman estimator for location or invariant Bayes estimators under right Haar measure---and hypothesis tests that are invariant with respect to nuisance group actions. In those classical works, the statistical model is typically well-specified and the group acts transitively on the parameter space to reduce risk. Our perspective differs in two important respects: first, we study the deterministic optimization landscape of a divergence functional under model misspecification ($p \notin \cQ_G$); second, the symmetry group $H$ of the standardized distribution and the transport group $G$ indexing the variational family need not coincide, and their interaction via conjugation governs which parameters are unbiasedly estimated. In a complementary direction, geometric deep learning and equivariant neural networks \citep{cohen2016group, weiler2019general, bronstein2021geometric} build group symmetries directly into function approximators and normalizing flows \citep{kohler2020equivariant}; rather than constraining the variational family to produce strictly $H$-invariant distributions (which would preclude estimating symmetry-breaking parameters $g^*$), our work uses $H$-invariant base distributions pushed forward by $G$ to estimate the unknown transformation $g^*$.

\textbf{Information geometry and optimization on Lie groups.} Differential-geometric approaches to statistical estimation have been extensively studied in information geometry \citep{murray1993differential, amari2016information}, where the Fisher--Rao metric defines a natural Riemannian structure on parametric families, and in optimization on smooth manifolds and Lie groups \citep{absil2008optimization, boumal2023introduction}. When a variational family is indexed by a Lie group $G$, left-trivialization identifies the tangent space $T_g G$ at any parameter $g$ with the Lie algebra $\gA = T_e G$. By equipping $\gA$ with an $\Ad_H$-invariant inner product (guaranteed by the compactness of $H$), we connect Riemannian first-order stationarity directly to the decomposition of $\gA$ into irreducible representations of $H$.

\textbf{Symmetry in variational inference.} The most direct predecessors of the present work are \citet{margossian2024variational} and \citet{margossian2025generalized}, who proved that location-scale variational inference exactly recovers the target mean under even symmetry and the target correlation matrix under elliptical symmetry across general $f$-divergences, and \citet{zellinger2026even}, who extended even-symmetry location recovery to forward KL and $\alpha$-divergences without log-concavity assumptions. Complementing these positive results, \citet{margossian2025variational} established an impossibility theorem for factorized (mean-field) variational approximations of non-factorized targets, proving that a mean-field family can accurately recover at most one of three canonical dispersion summaries (marginal variances, marginal precisions, or generalized variance). This impossibility underscores the necessity of non-factorized group-transport families such as full $\GL(d)$ covariance models. Concurrently with our work, \citet{marks2026symmetry} developed a global theory of symmetry-induced statistic recovery in variational inference, characterizing when global minimizers inherit target symmetries and applying their framework to von Mises--Fisher distributions on the sphere, while \citet{ghafourpour2026variational} investigated radial transport families. Our differential-geometric approach is directly complementary: by linearizing the conjugation action of $H$ on $G$ via the adjoint representation, we provide an explicit subspace decomposition of the transport Lie algebra $\gA$ that simultaneously explains both why symmetry-breaking gradient directions vanish and why symmetry-preserving directions ($\gA_{\text{bias}}^H$) exhibit estimation bias unless calibrated by moment matching.

%%===========================================================================
\section{Problem description}
\label{sec:problem}
%% ===========================================================================

In parametric variational inference \citep{wainwright2008graphical, ranganath2014black, blei2017variational}, one approximates an intractable target distribution $p \in \cP(\cX)$ on a sample space $\cX$ by selecting a candidate distribution $q_\theta$ from a parameterized family $\cQ_\Theta = \{q_\theta \in \cP(\cX) : \theta \in \Theta\}$ that minimizes a divergence functional over $\theta \in \Theta$. In this paper, we study structured variational families generated by the smooth action of a \emph{transport} Lie group $G$ on $\cX$, in the presence of a \emph{symmetry} group $H$ acting isometrically on $\cX$. For any $g \in G$ (respectively, $h \in H$) and any point $x \in \cX$, we denote the group action on the sample space by $g \cdot x \in \cX$ (respectively, $h \cdot x \in \cX$).

Viewing both $G$ and $H$ through their actions on the common domain $\cX$ serves two fundamental purposes. First, any measurable action on $\cX$ canonically induces a pushforward action on probability measures over $\cX$: transforming a random variable $Z \sim q_0$ via $X = g \cdot Z$ transports its law to $q_g = g_\sharp q_0$. Second, because both $G$ and $H$ act as transformations of $\cX$, they can be composed within the ambient group $\mathrm{Aut}(\cX)$ of diffeomorphisms of $\cX$, even when $G$ and $H$ are defined a priori as distinct algebraic groups. As we show throughout the paper, the non-commutativity of these two actions on $\cX$---captured algebraically by the conjugation of transports in $G$ by symmetries in $H$---is the central mechanism governing parameter identifiability and unbiased estimation in variational inference.

\subsection{Group transport variational inference}

We begin by formalizing variational families generated by pushing forward a fixed reference measure $q_0$ under a Lie group action. From a statistical standpoint, this construction generalizes classical location-scale transformation models \citep{eaton1989group} to arbitrary smooth manifolds and Lie groups; from a machine learning standpoint, it corresponds to a normalizing flow \citep{rezende2015variational, papamakarios2021normalizing} whose transformation maps are structured as a Lie group $G$. Throughout the paper, we restrict attention to distributions admitting strictly positive, smooth (or piecewise smooth, locally Lipschitz) densities with respect to the canonical reference measure on $\cX$ (Lebesgue measure when $\cX = \RR^d$, or the Riemannian volume measure when $\cX$ is a Riemannian manifold), and we identify each probability measure with its density.

\begin{definition}[Transport variational inference] \label{def:transport-vi}
Let $\cX$ be a smooth \emph{data domain} (such as Euclidean space $\RR^d$ or a smooth Riemannian manifold) and let $G$ be a Lie group acting smoothly and effectively on $\cX$ by diffeomorphisms, written $(g, x) \mapsto g \cdot x \in \cX$ (where acting \emph{effectively} simply means that distinct group elements $g \neq g'$ define distinct transformations on $\cX$). We call $G$ the \emph{parameter transport group}. Given a base distribution $q_0 \in \cP(\cX)$ and a target data distribution $p \in \cP(\cX)$, both with positive densities on $\cX$, we define:
\begin{itemize}[leftmargin=*]
  \item[-] The \emph{transported distribution} $q_g = g_\sharp q_0$: the law of $X = g \cdot Z$ when $Z \sim q_0$, whose density is given by the change-of-variables formula
  \[
  q_g(x) = q_0(g^{-1} \cdot x) \, |\det J_{g^{-1}}(x)|\, ,
  \]
  where $J_{g^{-1}}(x)$ denotes the Jacobian of the inverse map $x \mapsto g^{-1} \cdot x$ with respect to the reference volume measure on $\cX$ (in particular, $|\det J_{g^{-1}}(x)| = 1$ whenever $G$ acts by isometries).
  \item[-] The \emph{transport variational family} $\cQ_G = \{q_g \in \cP(\cX) \mid g \in G\}$.
  \item[-] The \emph{variational inference objective} $\cL : G \to \RR$, defined via the Kullback--Leibler divergence:
  \[
  \mathcal{L}(g) = \KL(q_g \| p) = \E_{x \sim q_g} \left[ \log q_g(x) - \log p(x) \right]\, .
  \]
\end{itemize}
\end{definition}

\paragraph{Illustrative examples of transport families.}
To ground Definition~\ref{def:transport-vi}, consider three classical parametric families obtained as group pushforwards of a standard base distribution $q_0$:
\begin{itemize}[leftmargin=*]
  \item[-] \textbf{Gaussian location family}: On $\cX = \RR^d$, let $G = (\RR^d, +)$ act by translation $g \cdot x = x + \mu_g$. Taking the standard Gaussian base $q_0 = \cN(0, I_d)$ yields the homoscedastic location family
  \[
  \cQ^{\text{mean}}_G = \{\cN(\mu_g, I_d) \mid \mu_g \in \RR^d \}\, .
  \]
  \item[-] \textbf{Gaussian scale/covariance family}: On $\cX = \RR^d$, let the general linear group $G = \GL(d)$ act linearly by $g \cdot x = A_g x$. If $Z \sim q_0 = \cN(0, I_d)$, then $A_g Z \sim \cN(0, A_g A_g^\top)$, generating the centered Gaussian covariance family
  \[
  \cQ^{\text{cov}}_G = \{\cN(0, A_g A_g^\top) \mid A_g \in \GL(d) \}\, .
  \]
  \item[-] \textbf{Directional von Mises--Fisher family on the sphere}: On the unit sphere $\cX = \cS^2 \subset \RR^3$, let $G = \SO(3)$ act by rotations $g \cdot x = R_g x$. Taking a von Mises--Fisher base distribution $q_0(x) = c_\kappa \exp(\kappa e_3^\top x)$ ($\kappa > 0$) concentrated around the polar axis $e_3 = (0,0,1)^\top$ (the North pole) yields the rotational family
  \[
  \cQ^{\text{ax}}_G = \{c_\kappa \exp(\kappa v_g^\top x) \mid v_g = R_g e_3,\, R_g \in \SO(3) \}\, ,
  \]
  indexed by the transported modal direction $v_g = R_g e_3 \in \cS^2$.
\end{itemize}

When the target distribution $p$ is itself generated by transporting a centered distribution $\tilde p$ via an unknown true transformation $g^* \in G$, a fundamental statistical question is whether minimizing $\cL(g)$ recovers $g^*$---or, at a first-order differential level, whether $g^*$ is a stationary point of $\cL$. Because $G$ is a smooth manifold rather than a flat vector space in general, first-order stationarity is characterized through the \emph{Riemannian gradient} $\nabla_G \cL(g^*)$, which is an element of the Lie algebra $\fG = T_e G$ tangent to $G$ at the identity (defined formally in Section~\ref{sec:lie}). By left-trivialization, $\nabla_G \cL(g)$ governs the first-order variation of $\cL$ along one-parameter subgroup perturbations $t \mapsto g \exp(t u)$ for $u \in \fG$:
\[
\cL\big(g \cdot \exp(tu)\big) = \cL(g) + t\langle\nabla_G\cL(g), u\rangle_\fG + o(t)\, .
\]
Thus, $g^*$ is a first-order critical point of the variational loss if and only if $\nabla_G \cL(g^*) = 0$ in $\fG$.

\textbf{Remarks on pushforwards, overparametrization, and notation.}
\begin{itemize}[leftmargin=*]
  \item[-] \emph{Pushforward versus pullback}: In measure-theoretic terminology, the transported distribution $q_g = g_\sharp q_0$ is the \emph{pushforward} of $q_0$ by the action of $g$. Conversely, standardizing a distribution $p$ by applying the inverse transformation $(g^*)^{-1}$ defines its \emph{pullback} $\tilde p = (g^*)^{-1}_\sharp p$ (Definition~\ref{def:symmetry}).
  \item[-] \emph{Overparametrization and coset identifiability}: In several natural families---such as the Gaussian covariance family ($G = \GL(d)$) and the spherical axis family ($G = \SO(3)$)---the map $g \mapsto q_g$ is not injective: distinct group elements $g \neq g'$ can produce the exact same distribution $q_g = q_{g'}$ (for instance, whenever $A_g A_g^\top = A_{g'} A_{g'}^\top$, which occurs when $A_{g'} = A_g U$ for $U \in \Or(d)$). In group-theoretic terminology, the set of transformations $n \in G$ that leave the base distribution unchanged ($q_n = q_0$) forms a subgroup called the \emph{stabilizer} of $q_0$, and the equivalence class of all parameters $\{g n : n \text{ in the stabilizer}\}$ that generate the same distribution $q_g$ is called a \emph{left coset}. In such settings, $g^*$ is identifiable only up to this stabilizer subgroup, while the induced distribution $q_{g^*}$ and the coset parameter on the quotient space (Section~\ref{sec:identifiability}) remain well-defined.
  \item[-] \emph{Notation}: Where no ambiguity arises, we identify a group element $g \in G$ directly with its concrete transformation parameter, writing $\mu \in \RR^d$ for translations, $A \in \GL(d)$ for invertible matrices, or $R \in \SO(d)$ for rotations.
\end{itemize}

\subsection{Distribution symmetries}

To understand when variational inference can recover $g^*$ despite functional misspecification between $q_0$ and the target distribution, we formalize the notion of distributional invariance under a group $H$ of isometries of $\cX$, and define \emph{conjugate symmetry} around a transport parameter $g^* \in G$. Statistically, conjugate symmetry expresses the familiar modeling assumption that a random variable $X \sim p$ becomes symmetric around the origin after being standardized by $(g^*)^{-1}$, i.e., the standardized variable $\xi = (g^*)^{-1} \cdot X \sim \tilde p$ has an $H$-invariant distribution.

\begin{definition}[$H$-symmetry and conjugate symmetry] \label{def:symmetry}
Let $H$ be a non-trivial compact group acting effectively on $\cX$ by isometries (written linearly as $h \cdot x = M_h x$ with $M_h \in \Or(d)$ when $\cX \subseteq \RR^d$).
\begin{itemize}[leftmargin=*]
  \item[-] A distribution $q \in \cP(\cX)$ is \emph{$H$-symmetric} (or \emph{$H$-invariant}) if $q(h\cdot z) = q(z)$ for all $h \in H$ and $z \in \cX$.
  \item[-] Given a target distribution $p \in \cP(\cX)$ and a transport element $g^* \in G$, the \emph{pullback distribution} $\tilde p = (g^*)^{-1}_\sharp p$ is the law of the standardized variable $\xi = (g^*)^{-1} \cdot X$ when $X \sim p$, with density
  \[
  \tilde p(\xi) = p(g^* \cdot \xi) \, |\det J_{g^*}(\xi)|\, .
  \]
  \item[-] A distribution $p \in \cP(\cX)$ has \emph{$H$-conjugate symmetry around $g^* \in G$} if its pullback $\tilde p$ under $g^*$ is $H$-symmetric, meaning that for all $h \in H$ and $\xi \in \cX$,
  \[
  \tilde p(h \cdot \xi) = \tilde p(\xi)\, .
  \]
\end{itemize}
\end{definition}

By construction, the pullback operation inverts the transport pushforward: if $X \sim p$, then $(g^*)^{-1} \cdot X \sim \tilde p$, and conversely if $\xi \sim \tilde p$, then $g^* \cdot \xi \sim p$. Thus, $p$ has $H$-conjugate symmetry around $g^*$ precisely when $p$ is obtained by transporting an $H$-symmetric reference distribution $\tilde p$ by $g^*$.

\paragraph{Examples of $H$-symmetry and conjugate symmetry.}
\begin{itemize}[leftmargin=*]
  \item[-] \emph{Central reflection symmetry around a location parameter}: Let $\cX = \RR^d$, let $H = \Z/2\Z$ act by central reflection $\{I_d, -I_d\}$, and let $G = (\R^d, +)$ act by translation $x \mapsto x + \mu$. A distribution $q \in \cP(\R^d)$ is $H$-symmetric if and only if it is even, $q(-z) = q(z)$. For a target distribution $p \in \cP(\R^d)$ and a location parameter $\mu^* \in \R^d$, the pullback is the centered density $\tilde p(\xi) = p(\mu^* + \xi)$, corresponding to the centered random variable $\xi = X - \mu^*$. Therefore, $p$ has $H$-conjugate symmetry around $\mu^*$ if and only if $\tilde p(\xi) = \tilde p(-\xi)$, or equivalently $p(\mu^* + \xi) = p(\mu^* - \xi)$ for all $\xi \in \R^d$. This recovers the classical setting of location estimation for centrally symmetric distributions studied in \citet{margossian2024variational}.
  \item[-] \emph{Canonical isometry groups in statistics}: Beyond central reflection $\Z/2\Z$, several compact isometry groups $H$ arise naturally in statistical modeling: the full orthogonal group $\Or(d)$ on $\RR^d$ (characterizing spherical and radial distributions, or elliptical distributions after affine standardization); the symmetric group $\Sn$ acting on $\RR^n$ by coordinate permutation (characterizing exchangeable random vectors); the dihedral groups $D_k$ on $\RR^2$ (characterizing discrete rotational and reflectional symmetries, such as $\ell_m$-norm contours); and the planar rotation group $\SO(2) \subset \SO(3)$ fixing a polar axis $e_3$ on the sphere $\cS^2$ (characterizing axially symmetric directional distributions).
\end{itemize}

\textbf{Remark (Geometric meaning of conjugate symmetry).} Written directly in terms of the target density $p$ without reference to the pullback $\tilde p$, $H$-conjugate symmetry around $g^*$ states that for every $h \in H$ and $x \in \cX$,
\[
p\big((g^* \circ h \circ (g^*)^{-1}) \cdot x\big) \, \big|\det J_{g^* \circ h \circ (g^*)^{-1}}(x)\big| = p(x)\, .
\]
That is, while the standardized distribution $\tilde p$ is invariant under the base symmetry group $H$ centered at the origin, transporting $\tilde p$ by $g^*$ shifts its symmetries to the \emph{conjugated subgroup} $g^* H (g^*)^{-1}$. Consequently, specifying a symmetric variational inference problem involves an inherently coupled modeling choice between the symmetry group $H$ of the standardized data and the transport group $G$ used to relocate the distribution. We formalize this coupling next.

\subsection{Transport variational inference with symmetry}

We now combine the transport family of Definition~\ref{def:transport-vi} with the symmetry structure of Definition~\ref{def:symmetry}. Suppose that the variational base distribution $q_0$ is $H$-symmetric, and that the target distribution $p$ has $H$-conjugate symmetry around an unknown true parameter $g^* \in G$. Even when the functional form of $q_0$ differs substantially from that of the pullback $\tilde p$ (for instance, when $q_0$ is Gaussian while $\tilde p$ is heavy-tailed or multimodal), both distributions share the same invariance group $H$. Our central question is whether this shared symmetry suffices for variational inference to identify $g^*$, and specifically whether $\nabla_G \cL(g^*) = 0$.

\citet{margossian2024variational} analyzed two special instances of this problem (corresponding to Examples~1 and~4 below), showing that the variational gradient vanishes at the true location parameter under reflection symmetry, but exhibits a nonzero radial scaling bias when estimating a covariance matrix under elliptical symmetry. Our goal is to establish a general theory covering arbitrary pairs $(G, H)$ that explains both phenomena and characterizes $\nabla_G \cL(g^*)$ completely.

\begin{definition}[Transport variational inference with symmetry] \label{def:vi-symmetry}
In the setting of Definition~\ref{def:transport-vi}, let $H$ be a non-trivial compact group acting effectively on $\cX$ by isometries (the \emph{symmetry group}). Assume that:
\begin{enumerate}[label=(\roman*), leftmargin=*]
  \item the variational base distribution $q_0 \in \cP(\cX)$ is $H$-symmetric (with stabilizer $\mathrm{Stab}_G(q_0) = G \cap H$ in $G$);
  \item the target distribution $p \in \cP(\cX)$ has $H$-conjugate symmetry around a true parameter $g^* \in G$; and
  \item the transport group $G$ is \emph{normalized} by $H$ within $\mathrm{Aut}(\cX)$, meaning that $\Phi_h(g) \coloneqq h \circ g \circ h^{-1} \in G$ for all $g \in G$ and $h \in H$---in plain terms, expressing any parameter transformation $g \in G$ in a coordinate system rotated or reflected by $h \in H$ yields another valid parameter transformation inside $G$.
\end{enumerate}
Symmetric transport variational inference then consists in solving the optimization problem
\[
\min_{g\in G}\cL(g) = \min_{g \in G}\KL(q_g \| p)\, .
\]
\end{definition}

\textbf{Remark (Why $G$ must be normalized by $H$).} Condition~(iii) in Definition~\ref{def:vi-symmetry} requires that conjugating any transport $g \in G$ by any symmetry $h \in H$ produces a transformation $\Phi_h(g) = h \circ g \circ h^{-1}$ that still belongs to the variational transport group $G$. This structural consistency condition has a clear geometric and statistical motivation.

Consider perturbing the variational parameter away from $g^*$ by a relative transport $g \in G$, evaluating $\cL(g^* g) = \KL(q_g \| \tilde p)$. Because both $q_0$ and $\tilde p$ are $H$-invariant, rotating the perturbation by any symmetry $h \in H$ produces a symmetrically equivalent perturbation $h \circ g \circ h^{-1}$ that achieves the exact same KL divergence against $\tilde p$. However, for this symmetry of the objective to constrain the gradient within our parameter space $G$, the rotated perturbation $h \circ g \circ h^{-1}$ must actually lie in $G$. If $h \circ g \circ h^{-1} \notin G$, the variational family truncates the symmetry orbit of the perturbation. For example, if $\tilde p$ is rotationally invariant in the plane ($H = \SO(2)$) but the variational family only allows horizontal translations along the $x_1$-axis, rotating a horizontal translation by $90^\circ$ produces a vertical translation along the $x_2$-axis, which lies outside $G$. Requiring $G$ to be normalized by $H$ ensures that the variational parameter space is closed under the symmetry transformations we seek to exploit, and guarantees that the quotient parameter space $\cM = G / (G \cap H)$ (Section~\ref{sec:identifiability}) and the adjoint representation on the Lie algebra $\gA$ (Section~\ref{sec:lie}) are well-defined.

\paragraph{Guiding examples.}
Throughout the paper, we trace six concrete examples that span location estimation, directional statistics, covariance estimation, non-Gaussian orientation estimation, and exchangeable models. These examples allow every algebraic object in our theory to be interpreted directly in familiar statistical terms:

\begin{enumerate}[leftmargin=*]
  \item \textbf{Reflection (Univariate symmetric location model)}: On the real line $\cX = \RR$, let $H = \Z / 2\Z = \{1, -1\}$ act by sign reflection $x \mapsto \pm x$, and let $G = (\RR, +)$ act by translation $x \mapsto x + \mu$. Choosing the standard normal base $q_0 = \cN(0, 1)$ generates the unit-variance Gaussian location family $\cQ_\R = \{\cN(\mu, 1) \mid \mu \in \R\}$. Any target distribution $p$ that is symmetric around a true location $\mu^* \in \R$---such as a symmetric bimodal Gaussian mixture $p = \tfrac{1}{2}\cN(\mu^* + \delta, 1) + \tfrac{1}{2}\cN(\mu^* - \delta, 1)$, a Student-$t$ distribution centered at $\mu^*$, or a Laplace distribution---has $H$-conjugate symmetry around $\mu^*$ even though $p \notin \cQ_\R$.

  \item \textbf{Euclidean (Multivariate radially symmetric location model)}: On $\cX = \RR^d$ ($d \ge 2$), let $H = \Or(d)$ be the orthogonal group of rotations and reflections around the origin, and let $G = (\RR^d, +)$ act by translation $x \mapsto x + \mu$. With $q_0 = \cN(0, I_d)$, the variational family is $\cQ_{\R^d} = \{\cN(\mu, I_d) \mid \mu \in \R^d\}$, while the target $p(x) = \tilde p(x - \mu^*)$ can be any radially symmetric distribution $\tilde p(\xi) = f(\|\xi\|_2)$ shifted to center $\mu^* \in \RR^d$.

  \item \textbf{Sphere (Directional estimation for axially symmetric distributions)}: On the unit sphere $\cX = \cS^2 \subset \R^3$, let $H = \SO(2)$ be the subgroup of $\SO(3)$ consisting of rotations around the polar axis $e_3 = (0,0,1)^\top$ (the North pole), and let $G = \SO(3)$ act by 3D rotations $x \mapsto Rx$. Using a von Mises--Fisher base distribution $q_0(x) = c_\kappa \exp(\kappa e_3^\top x)$ ($\kappa > 0$) produces the variational family $\cQ_{\SO(3)} = \{q_R(x) = c_\kappa \exp(\kappa (R e_3)^\top x) \mid R \in \SO(3)\}$, which is effectively indexed by the modal direction $v = Re_3 \in \cS^2$. The target $p(x) = \tilde p((R^*)^\top x)$ is obtained by rotating any rotationally symmetric spherical distribution $\tilde p$ around $e_3$ so that its axis of symmetry aligns with $v^* = R^* e_3$.

  \item \textbf{Covariance (Elliptical scatter and covariance estimation)}: On $\cX = \RR^d$ ($d \ge 2$), let $H = \Or(d)$ act by orthogonal transformations and let $G = \GL(d)$ act linearly by $x \mapsto Ax$. Taking $q_0 = \cN(0, I_d)$ yields the zero-mean Gaussian family $\cQ^{\text{cov}}_{\GL(d)} = \{\cN(0, A A^\top) \mid A \in \GL(d)\}$, which parameterizes all symmetric positive-definite covariance matrices $\Sigma = A A^\top \in \mathrm{SPD}(d)$. If the target $p(x) = |\det A^*|^{-1} \tilde p((A^*)^{-1} x)$ is generated by applying an invertible linear transformation $A^*$ to any spherically symmetric distribution $\tilde p$ (such as a multivariate $t$ or Kotz-type elliptical distribution), then $p$ is an elliptical distribution with scatter structure $\Sigma^* = A^* (A^*)^\top$ and possesses $\Or(d)$-conjugate symmetry around $A^*$.

  \item \textbf{Square (Orientation estimation under discrete dihedral symmetry)}: On the plane $\cX = \RR^2$, let $H = D_4$ be the dihedral symmetry group of the square (generated by $90^\circ$ rotations and coordinate-axis reflections), and let $G = \SO(2)$ act by continuous planar rotations $R_\theta$. A natural $D_4$-invariant base distribution is the generalized Gaussian $q_0(x) \propto \exp(-\gamma_m \|x\|_m^m)$ with $m > 2$ (for example, $m = 4$), whose level sets are $\ell_m$-balls with fourfold corners. Rotating $q_0$ by $R_\theta \in \SO(2)$ generates the orientation family $\cQ^{\text{square}}_{\SO(2)} = \{q_0(R_\theta^\top x) \mid R_\theta \in \SO(2)\}$ modeling distributions whose square-symmetric contours are tilted by an angle $\theta$, against a target $p(x) = \tilde p(R_{\theta^*}^\top x)$ tilted by a true angle $\theta^*$.

  \item \textbf{Permutation (Location estimation for exchangeable coordinates)}: On $\cX = \RR^n$, let $H = \Sn$ be the symmetric group of coordinate permutations $(x_1, \dots, x_n) \mapsto (x_{\sigma(1)}, \dots, x_{\sigma(n)})$, which characterizes exchangeable random variables, and let $G = (\RR^n, +)$ act by coordinate-wise translation $x \mapsto x + \mu$. For instance, taking an exchangeable product Laplace base $q_0(x) = c_\kappa \exp(-\kappa \|x\|_1)$ yields the location family $\cQ_{\R^n}^{\text{perm}} = \{c_\kappa\exp(-\kappa\|x - \mu\|_1) \mid \mu \in \R^n\}$, while the target $p(x) = \tilde p(x - \mu^*)$ is obtained by shifting any exchangeable distribution $\tilde p$ (potentially with dependent coordinates, such as an equicorrelated normal or exchangeable mixture) by a location vector $\mu^* \in \RR^n$.
\end{enumerate}

Table~\ref{tab:examples-preview} summarizes the sample spaces, transport groups, and symmetry groups for these six guiding examples.

\begin{table}[ht!]
\centering
\caption{Summary of sample spaces, transport groups, and symmetry groups for the six guiding examples.}
\label{tab:examples-preview}
{\small
\begin{tabular}{l|c|c|c|l}
\textbf{Example} & \begin{tabular}[c]{@{}c@{}}\textbf{Data}\\ \textbf{Domain} ($\cX$)\end{tabular} & \begin{tabular}[c]{@{}c@{}}\textbf{Transport}\\ \textbf{Group} ($G$)\end{tabular} & \begin{tabular}[c]{@{}c@{}}\textbf{Symmetry}\\ \textbf{Group} ($H$)\end{tabular} & \textbf{Symmetry Action} \\ \hline
1. Reflection & $\RR$ & $(\RR, +)$ & $\Z/2\Z$ & Sign flip ($x \mapsto -x$) \\
2. Euclidean & $\RR^d$ & $(\RR^d, +)$ & $\Or(d)$ & Rotations / reflections \\
3. Sphere & $\cS^2 \subset \RR^3$ & $\SO(3)$ & $\SO(2)$ & Axial rotations around $e_3$ \\
4. Covariance & $\RR^d$ & $\GL(d)$ & $\Or(d)$ & Orthogonal basis change \\
5. Square & $\RR^2$ & $\SO(2)$ & $D_4$ & $90^\circ$ rot.\ \& axis flips \\
6. Permutation & $\RR^n$ & $(\RR^n, +)$ & $\Sn$ & Coordinate permutations \\
\end{tabular}
}
\end{table}

%% ===========================================================================
\section{Group interactions and identifiability} \label{sec:identifiability}
%% ===========================================================================

Before analyzing first-order stationarity of the variational loss $\cL(g) = \KL(q_g \| p)$, we must address a foundational statistical question: \emph{which parameters are identifiable from the variational family $\cQ_G = \{q_g : g \in G\}$?} Because the base distribution $q_0$ is invariant under the symmetry group $H$, any transport $n \in G$ that also belongs to $H$ leaves $q_0$ unchanged ($q_n = q_0$), and more generally right-multiplying any $g \in G$ by $n \in G \cap H$ produces the identical distribution $q_{gn} = q_g$. Furthermore, even when $h \in H$ does not belong to $G$, composing transports in $G$ with symmetries in $H$ on the sample space $\cX$ generates a larger ambient group of transformations. In this section, we formalize this ambient group $K = GH$, analyze the conjugation action $\Phi_h(g) = h \circ g \circ h^{-1}$ across our six guiding examples, and classify the resulting statistical identifiability regimes on the quotient manifold $\cM = G / (G \cap H)$.

\subsection{Total transformation group and semidirect product structure}

When $G$ and $H$ both act on the sample space $\cX$, their compositions generate a combined group of transformations inside $\mathrm{Aut}(\cX)$. The algebraic structure of this total group is governed by the normalization condition $\Phi_h(g) = h \circ g \circ h^{-1} \in G$ introduced in Definition~\ref{def:vi-symmetry}.

\begin{definition}[Total transformation group] \label{def:total-group}
Let $K = GH \subseteq \mathrm{Aut}(\cX)$ denote the total group of transformations on $\cX$ generated by composing actions of $G$ and $H$. Every element $k \in K$ can be written as $k = g \circ h$ for some $g \in G$ and $h \in H$, acting on $x \in \cX$ by applying the symmetry $h$ first and the transport $g$ second:
\[
k \cdot x = g \cdot (h \cdot x)\, .
\]
Because $G$ is normalized by $H$, the composition of any two elements $k_2 = g_2 \circ h_2$ and $k_1 = g_1 \circ h_1$ in $K$ remains in $GH$ and obeys the semidirect multiplication law
\[
(g_2 \circ h_2) \circ (g_1 \circ h_1) = \big(g_2 \circ \Phi_{h_2}(g_1)\big) \circ (h_2 \circ h_1)\, ,
\]
where $\Phi_h(g) = h \circ g \circ h^{-1} \in G$ denotes the conjugation of $g$ by $h$. Depending on the intersection subgroup $N = G \cap H$, the total group $K$ takes one of three algebraic forms:
\begin{enumerate}[label=(\roman*), leftmargin=*]
  \item When $N = G \cap H = \{e\}$ (disjoint actions), $K$ is isomorphic to the internal semidirect product $G \rtimes_\Phi H$: every transformation $k \in K$ decomposes uniquely into a symmetry $h \in H$ followed by a transport $g \in G$, and composing two such transformations twists the first transport parameter $g_1$ by $h_2$ via $\Phi_{h_2}(g_1)$ (just as rotating an affine-shifted variable rotates its translation vector). For location families where $G = (\RR^d, +)$ acts by translation and $H$ acts linearly (Examples~1, 2, and~6), $K = \RR^d \rtimes H$ is the group of affine isometries generated by $H$ and translations (in particular, the full Euclidean group $\mathrm{E}(d) = \RR^d \rtimes \Or(d)$ in Example~2).
  \item When $H \subseteq G$, so that $N = H$ (internal symmetries), every symmetry is already a transport transformation, and the total group collapses to $K = G$ (as in Example~3 where $K = \SO(3)$ and Example~4 where $K = \GL(d)$).
  \item In the general overlapping case where $\{e\} \subsetneq N \subsetneq H$, $K$ is isomorphic to the quotient $(G \rtimes_\Phi H) / N$, where $N$ is embedded anti-diagonally via $n \mapsto (n, n^{-1})$ (as in Example~5, where $G = \SO(2)$ and $H = D_4$ intersect in $N = C_4$ and together generate the full planar orthogonal group $K = \Or(2)$).
\end{enumerate}
\end{definition}

\textbf{Remark (Closure under conjugation and well-posedness of $K$).} Definition~\ref{def:total-group} underscores why the normalization requirement $\Phi_h(g) \in G$ from Section~\ref{sec:problem} is algebraically indispensable. Rewriting the composition $h_2 \circ g_1$ as $(h_2 \circ g_1 \circ h_2^{-1}) \circ h_2 = \Phi_{h_2}(g_1) \circ h_2$ allows any sequence of alternating transports and symmetries to be re-ordered into a single symmetry followed by a single transport in $G$. If $G$ were not normalized by $H$, composing transports with symmetries would generate transformations outside $GH$, and differentiating the conjugation map at the identity (Section~\ref{sec:lie}) would produce tangent vectors pointing outside the Lie algebra $\gA$ of the variational family.

\subsection{The conjugation action in detail} \label{sec:conjugation-detail}

The conjugation map $\Phi : H \to \mathrm{Aut}(G)$, $h \mapsto \Phi_h$, where $\Phi_h(g) = h \circ g \circ h^{-1}$, is the fundamental algebraic bridge connecting the symmetry group $H$ to the transport group $G$. Before passing to the tangent space in Section~\ref{sec:lie}, we establish the basic structural properties of $\Phi$, examine its geometric meaning on $\cX$, and compute its action, kernel, and fixed-point subgroup across all six guiding examples.

\begin{proposition}[Conjugation as a group homomorphism] \label{prop:conjugation-homomorphism}
The map $\Phi: H \to \mathrm{Aut}(G)$, $h \mapsto \Phi_h$, is a group homomorphism: $\Phi_{h_1 h_2} = \Phi_{h_1} \circ \Phi_{h_2}$ for all $h_1, h_2 \in H$. Moreover, for each $h \in H$, the map $\Phi_h : G \to G$ is a Lie group automorphism, satisfying $\Phi_h(g_1 g_2) = \Phi_h(g_1)\Phi_h(g_2)$ and $\Phi_h(e_G) = e_G$.
\end{proposition}

\begin{proof}
For any $h_1, h_2 \in H$ and $g \in G$, associativity of composition in $\mathrm{Aut}(\cX)$ yields
\[
\Phi_{h_1 h_2}(g) = (h_1 h_2) \circ g \circ (h_1 h_2)^{-1} = h_1 \circ (h_2 \circ g \circ h_2^{-1}) \circ h_1^{-1} = \Phi_{h_1}(\Phi_{h_2}(g))\, .
\]
Similarly, for any $g_1, g_2 \in G$, inserting $h^{-1} \circ h = \mathrm{id}_{\cX}$ gives
\[
\Phi_h(g_1 g_2) = h \circ g_1 \circ (h^{-1} \circ h) \circ g_2 \circ h^{-1} = \Phi_h(g_1)\Phi_h(g_2)\, ,
\]
and $\Phi_h(e_G) = h \circ e_G \circ h^{-1} = e_G$. Smoothness of $\Phi_h$ follows from the smooth group structure of $K$.
\end{proof}

Two subgroups associated with the conjugation action $\Phi$ play a recurring role in our analysis:
\begin{itemize}[leftmargin=*]
  \item[-] The \emph{kernel of conjugation} (or \emph{centralizer} of $G$ in $H$),
  \[
  \ker(\Phi) = \{h \in H : \Phi_h(g) = g \;\; \forall g \in G\} = C_H(G)\, ,
  \]
  which consists of all symmetries in $H$ that commute with every transport in $G$ (meaning that transforming the coordinate system by $h$ leaves every parameter $g \in G$ unchanged). When $\ker(\Phi) = \{e_H\}$, the conjugation action of $H$ on $G$ is faithful: every non-trivial symmetry induces a non-trivial reparametrization of the transport group.
  \item[-] The \emph{fixed-point subgroup} (or \emph{$H$-centralizer} in $G$),
  \[
  G^H = \{g \in G : \Phi_h(g) = g \;\; \forall h \in H\}\, ,
  \]
  which consists of all transports in $G$ that commute with every symmetry in $H$. Because each $\Phi_h$ is a group automorphism fixing $e_G$, $G^H$ is a closed Lie subgroup of $G$. Statistically, $G^H$ consists of those parameter transformations (such as global isotropic scaling) that do not break the $H$-symmetry of the model: if $q_0$ is $H$-symmetric and $g \in G^H$, then $q_g$ remains $H$-symmetric.
\end{itemize}

\paragraph{Geometric interpretation.}
Conjugation describes how a change of reference frame by a symmetry $h \in H$ re-expresses a transport transformation $g \in G$. Suppose we wish to apply a transport $g$ in a coordinate system that has been transformed by $h$. Moving from the original coordinates to the $h$-transformed coordinates applies $h^{-1}$; performing the transport in those coordinates applies $g$; and returning to the original frame applies $h$. The resulting transformation on $x \in \cX$ is $(h \circ g \circ h^{-1}) \cdot x = \Phi_h(g) \cdot x$. Equivalently, for every $x \in \cX$, $g \in G$, and $h \in H$, the following diagram commutes:
\[
\begin{array}{ccc}
x & \xrightarrow{\quad g \quad} & g \cdot x \\
{\scriptstyle h}\downarrow \phantom{h} & & \phantom{h} \downarrow {\scriptstyle h} \\
h \cdot x & \xrightarrow{\;\Phi_h(g)\;} & h \cdot (g \cdot x) = \Phi_h(g) \cdot (h \cdot x)
\end{array}
\]
That is, transporting a point by $g$ and then applying the symmetry $h$ is identical to first applying the symmetry $h$ and then transporting by the conjugated parameter $\Phi_h(g) \in G$.

\paragraph{Detailed analysis of the guiding examples.}
We now compute the total group $K$, the conjugation map $\Phi_h$, the kernel $\ker(\Phi)$, and the fixed-point subgroup $G^H$ explicitly for each of our six guiding examples:

\begin{enumerate}[leftmargin=*]
\item \textbf{Reflection} ($G = (\RR, +)$, $H = \Z/2\Z = \{e, \sigma\}$ with $\sigma \cdot x = -x$).
The total transformation group is the one-dimensional semidirect product $K = \RR \rtimes (\Z/2\Z)$ of translations and sign reflections. Denoting translation by $T_\mu(x) = x + \mu$, conjugating $T_\mu$ by the reflection $\sigma$ yields
\[
(\sigma \circ T_\mu \circ \sigma^{-1})(x) = \sigma(T_\mu(-x)) = -(-x + \mu) = x - \mu = T_{-\mu}(x)\, .
\]
Thus $\Phi_\sigma(\mu) = -\mu$: reflecting the coordinate axis reverses the direction of the translation shift. The kernel $\ker(\Phi) = \{e\}$ is trivial, and the fixed-point subgroup is $G^H = \{\mu \in \RR : -\mu = \mu\} = \{0\}$, indicating that no nonzero translation preserves reflection symmetry around the origin.

\item \textbf{Euclidean} ($G = (\RR^d, +)$, $H = \Or(d)$ with $h \cdot x = Ux$ for $U \in \Or(d)$).
The total group is the full Euclidean group $K = \mathrm{E}(d) = \RR^d \rtimes \Or(d)$. Conjugating a translation $T_\mu(x) = x + \mu$ by an orthogonal transformation $U \in \Or(d)$ gives
\[
(U \circ T_\mu \circ U^{-1})(x) = U(U^{-1}x + \mu) = x + U\mu = T_{U\mu}(x)\, ,
\]
so $\Phi_U(\mu) = U\mu$: an orthogonal change of basis rotates or reflects the location parameter $\mu$ by the same matrix $U$. The conjugation orbits in $G = \RR^d$ are concentric spheres $\{\mu \in \RR^d : \|\mu\|_2 = r\}$. Because only the identity matrix fixes every vector in $\RR^d$, the action is faithful ($\ker(\Phi) = \{I_d\}$), and the fixed-point subgroup is trivial, $G^H = \{0\}$.

\item \textbf{Sphere} ($G = \SO(3)$, $H = \SO(2)$ embedded as planar rotations around $e_3$).
Because $H \subset G$, the total group collapses to $K = \SO(3)$. For an axial symmetry rotation $h = R_\alpha^{e_3} \in \SO(2)$ by angle $\alpha$ around $e_3$ and a general 3D rotation $g = R_\beta^v \in \SO(3)$ by angle $\beta$ around a unit axis $v \in \cS^2$, matrix conjugation yields
\[
\Phi_h(g) = R_\alpha^{e_3} R_\beta^v (R_\alpha^{e_3})^{-1} = R_\beta^{R_\alpha^{e_3} v}\, .
\]
Geometrically, conjugating a rotation by $R_\alpha^{e_3}$ preserves the rotation angle $\beta$ while rotating its axis $v$ around the North--South pole $e_3$ by angle $\alpha$. A rotation $R_\beta^v$ commutes with all rotations around $e_3$ if and only if its axis $v$ is parallel to $e_3$; hence the fixed-point subgroup is $G^H = \{R_\gamma^{e_3} : \gamma \in [0, 2\pi)\} = \SO(2) = H$, which coincides exactly with the internal symmetry group $N = G \cap H = H$.

\item \textbf{Covariance} ($G = \GL(d)$, $H = \Or(d)$ with $g \cdot x = Ax$).
Here again $H \subset G$, so the total group is $K = \GL(d)$. For an orthogonal matrix $U \in \Or(d)$ and an invertible matrix $A \in \GL(d)$, conjugation acts by orthogonal similarity:
\[
\Phi_U(A) = U A U^{-1} = U A U^\top\, .
\]
Statistically, conjugating $A$ by $U$ corresponds to rotating the coordinate axes of the covariance matrix $\Sigma = A A^\top \mapsto U \Sigma U^\top$, which rotates the principal axes of the ellipsoid while preserving its eigenvalues (and hence its trace and determinant). Any matrix $A \in \GL(d)$ satisfying $U A U^\top = A$ for all $U \in \Or(d)$ must remain invariant under all orthogonal rotations and coordinate reflections, which forces $A$ to have equal scaling in every direction and no off-diagonal couplings (formally, by Schur's lemma applied to the irreducible action of $\Or(d)$ on $\RR^d$). Thus the fixed-point subgroup is the one-dimensional subgroup of isotropic dilations, $G^H = \{\lambda I_d : \lambda \in \RR \setminus \{0\}\}$, whereas the kernel of $\Phi$ is $\ker(\Phi) = \{\pm I_d\}$.

\item \textbf{Square} ($G = \SO(2)$, $H = D_4$).
The planar rotation group $G = \SO(2)$ and the dihedral group $H = D_4$ intersect in the fourfold rotational subgroup $N = C_4 = \{R_{k\pi/2} : k = 0,1,2,3\}$ and together generate the full orthogonal group $K = \Or(2)$. The eight elements of $D_4$ consist of four rotations $r_k = R_{k\pi/2}$ and four reflections $s_k$. Because $\SO(2)$ is abelian, all rotations $r_k \in C_4$ commute with every $R_\theta \in \SO(2)$, giving $\Phi_{r_k}(R_\theta) = R_\theta$. By contrast, any reflection $s_k \in D_4$ reverses planar orientation and therefore inverts the rotation angle:
\[
\Phi_{s_k}(R_\theta) = s_k R_\theta s_k^{-1} = R_{-\theta}\, .
\]
Consequently, the conjugation homomorphism $\Phi : D_4 \to \mathrm{Aut}(\SO(2))$ factors through $D_4 / C_4 \cong \Z/2\Z$ with kernel $\ker(\Phi) = C_4$, and its fixed-point set in $\SO(2)$ is the discrete two-element subgroup $G^H = \{R_0, R_\pi\} \cong \Z/2\Z$.

\item \textbf{Permutation} ($G = (\RR^n, +)$, $H = \Sn$ with $\sigma \cdot x = P_\sigma x$).
The total group is the semidirect product $K = \RR^n \rtimes \Sn$ of translations and coordinate permutations. Conjugating a translation $T_\mu(x) = x + \mu$ by a permutation matrix $P_\sigma$ gives
\[
(P_\sigma \circ T_\mu \circ P_\sigma^{-1})(x) = P_\sigma(P_\sigma^{-1}x + \mu) = x + P_\sigma\mu = T_{P_\sigma\mu}(x)\, ,
\]
so $\Phi_\sigma(\mu) = P_\sigma \mu$: coordinate permutations act on the location parameter by permuting its entries. The conjugation orbit of $\mu \in \RR^n$ consists of all vectors sharing the same empirical multiset of coordinates. The action is faithful ($\ker(\Phi) = \{e\}$), and a location vector $\mu$ is fixed by all permutations if and only if all of its coordinates are equal:
\[
G^H = \{\mu \in \RR^n : P_\sigma \mu = \mu \;\; \forall \sigma \in \Sn\} = \mathrm{span}(\ones)\, ,
\]
where $\ones = (1, \dots, 1)^\top \in \RR^n$ denotes the all-ones vector. Statistically, shifting an exchangeable distribution by $c\ones$ shifts the common location of all coordinates simultaneously while preserving exchangeability.
\end{enumerate}

\paragraph{Conjugation and distributional invariance of the variational objective.}
The following proposition explains why the conjugation map $\Phi_h$ directly governs the symmetries of the variational loss landscape $\cL(g) = \KL(q_g \| p)$ around $g^*$.

\begin{proposition}[Conjugation and transport of symmetric distributions] \label{prop:conjugation-distribution}
Assume that the base distribution $q_0 \in \cP(\cX)$ is $H$-symmetric, and let $g \in G$. For any $h \in H$, let $q_{hg} \coloneqq h_\sharp q_g$ denote the pushforward of $q_g$ under the symmetry $h$, i.e., the distribution of $h \cdot (g \cdot Z)$ when $Z \sim q_0$ (note that the notation $q_{hg}$ refers to the pushforward measure $(h \circ g)_\sharp q_0$, even when $h \notin G$). Then:
\begin{enumerate}[leftmargin=*]
\item[(i)] \emph{Pushforward equivariance}: $q_{\Phi_h(g)} = q_{hg}$, meaning that $q_{\Phi_h(g)}(x) = q_g(h^{-1} \cdot x)$ for all $x \in \cX$.
\item[(ii)] \emph{Loss under conjugated transport}: $\cL(\Phi_h(g)) = \KL(q_{hg} \| p)$ for all $h \in H$ and $g \in G$.
\item[(iii)] \emph{Conjugation invariance of the relative loss around $g^*$}: If, in addition, the target $p$ has $H$-conjugate symmetry around $g^* \in G$, then for every relative transport $g \in G$ and every symmetry $h \in H$,
\[
\cL\big(g^* \Phi_h(g)\big) = \cL(g^* g)\, .
\]
\end{enumerate}
\end{proposition}

\begin{proof}
To prove (i), observe that $q_{\Phi_h(g)}$ is the law of $\Phi_h(g) \cdot Z = (h \circ g \circ h^{-1}) \cdot Z$ when $Z \sim q_0$. Because $q_0$ is $H$-symmetric, the random variable $Z' = h^{-1} \cdot Z$ has the same distribution $q_0$ as $Z$. Therefore,
\[
\Phi_h(g) \cdot Z = h \cdot \big(g \cdot (h^{-1} \cdot Z)\big) \stackrel{d}{=} h \cdot (g \cdot Z) \sim q_{hg}\, ,
\]
which establishes $q_{\Phi_h(g)} = q_{hg}$, and part~(ii) follows immediately from the definition of $\cL$. For part~(iii), recall that $p = (g^*)_\sharp \tilde p$ where $\tilde p$ is $H$-symmetric ($\tilde p = h_\sharp \tilde p$). By the invariance of the KL divergence under the diffeomorphism $(g^*)^{-1}$, we have $\cL(g^* g') = \KL(q_{g^* g'} \| p) = \KL(q_{g'} \| \tilde p)$ for any $g' \in G$. Applying this identity with $g' = \Phi_h(g)$ and invoking part~(i) together with the $H$-invariance of $\tilde p$ yields
\[
\cL\big(g^* \Phi_h(g)\big) = \KL\big(q_{\Phi_h(g)} \,\big\|\, \tilde p\big) = \KL\big(h_\sharp q_g \,\big\|\, \tilde p\big) = \KL\big(h_\sharp q_g \,\big\|\, h_\sharp \tilde p\big) = \KL(q_g \| \tilde p) = \cL(g^* g)\, ,
\]
where the penultimate equality again uses the invariance of KL divergence under the diffeomorphism $h$.
\end{proof}

Proposition~\ref{prop:conjugation-distribution} provides the key geometric insight into how symmetry shapes the variational objective. Part~(i) shows that because $q_0$ is $H$-symmetric, applying a symmetry $h$ to the transported distribution $q_g$ is equivalent to replacing the parameter $g$ by its conjugate $\Phi_h(g) \in G$---at the distribution level, $q_g$ and $q_{\Phi_h(g)}$ have the exact same shape, viewed from an $h$-rotated angle. Part~(iii) then shows that when $p$ is $H$-conjugate symmetric around $g^*$, stepping away from $g^*$ by a relative transport $g$ or by its conjugated version $\Phi_h(g)$ yields the exact same KL divergence, so the local optimization landscape around $g^*$ cannot distinguish between conjugate relative transports and $\cL$ is constant along every conjugation orbit $\{g^* \Phi_h(g) : h \in H\}$. Notice that this conjugation invariance $\cL(g^* \Phi_h(g)) = \cL(g^* g)$ is far richer than the right-invariance $\cL(gn) = \cL(g)$ for $n \in N = G \cap H$: right-invariance under $N$ merely reflects parameter redundancy \emph{within} $G$ (since $q_{gn} = q_g$), whereas conjugation invariance captures how symmetries in $H$ act on $G$ \emph{from the outside} to relate distinct distributions $q_{g^* g}$ and $q_{g^* \Phi_h(g)}$ that are symmetric counterparts of one another around $g^*$.

\subsection{Quotient space and statistical identifiability regimes}

We now formalize the parameter redundancy induced by the overlap between the transport group $G$ and the symmetry group $H$. Because $q_0$ is $H$-invariant, any transformation in the total group $K = GH$ that differs on the right by an element of $H$ produces the exact same pushforward distribution on $\cX$. Restricting to the transport group $G$, two parameters $g, g' \in G$ satisfy $g' = gn$ for some $n \in N \coloneqq G \cap H$ if and only if they belong to the same left coset $gN$ in the quotient space
\[
\cM \coloneqq G / N \cong K / H\, , \qquad N = G \cap H\, .
\]
In statistical terms, the subgroup $\mathrm{Stab}_G(q_0) \coloneqq \{g \in G : g_\sharp q_0 = q_0\}$ (the \emph{stabilizer} of $q_0$ in $G$) is the set of redundant parameter transformations that leave the base distribution $q_0$ unchanged, and forming the \emph{quotient space} $\cM = G/N$ simply means grouping together all parameters $\{gn : n \in N\}$ that produce the exact same distribution $q_g$ into a single equivalence class $[g] = gN$ (just as a covariance matrix $\Sigma = A A^\top$ groups together all matrix square roots $\{A U : U \in \Or(d)\}$). Assuming that the stabilizer of $q_0$ within $G$ is precisely $N = G \cap H$ (i.e., $\mathrm{Stab}_G(q_0) = G \cap H$), the map $[g] = gN \mapsto q_g$ is injective on $\cM$, so the smooth coset manifold $\cM = G/N$ represents the space of \emph{statistically identifiable parameters} of the variational family $\cQ_G$. Depending on the structure of the intersection $N = G \cap H$, our examples fall into three distinct identifiability regimes:

\paragraph{1. Separable symmetry ($N = G \cap H = \{e\}$).}
When the transport group $G$ and the symmetry group $H$ intersect only at the identity transformation, the quotient manifold coincides with the transport group itself, $\cM \cong G$. In this regime, no non-trivial transport leaves $q_0$ invariant, so the parameter $g \in G$ is \textbf{globally identifiable} from $q_g$:
\begin{itemize}[leftmargin=*]
\item[-] \textbf{Examples~1, 2, and~6 (Reflection, Euclidean, and Permutation location families)}: In all three translation models, $G = (\RR^d, +)$ (or $\RR^n$) contains no non-trivial linear isometries around the origin, so $N = G \cap H = \{0\}$ and $\cM \cong \RR^d$ (or $\RR^n$). Every location parameter $\mu$ uniquely indexes a distinct distribution $q_\mu \in \cQ_G$. Crucially, however, global parameter identifiability on $\cM$ does not by itself guarantee that $\mu^*$ minimizes $\cL(\mu)$ under model misspecification: as we show in Section~\ref{sec:main-result}, Examples~1 and~2 are always stationary at $\mu^*$, whereas Example~6 can exhibit an additive location bias along the diagonal $\mathrm{span}(\ones)$.
\end{itemize}

\paragraph{2. Internal symmetry ($H \subset G$, so $N = H$).}
When the symmetry group $H$ is a subgroup of the transport group $G$, the total group collapses to $K = G$ and the identifiable parameter space is the homogeneous quotient manifold $\cM = G / H$ of left cosets $[g] = gH$. In this regime, the group element $g \in G$ is \textbf{overparametrized}, and only the coset $[g] \in G/H$ is statistically identifiable:
\begin{itemize}[leftmargin=*]
\item[-] \textbf{Example~3 (Sphere)}: With $G = \SO(3)$ and $H = N = \SO(2)$ (rotations around $e_3$), the quotient manifold is $\cM = \SO(3)/\SO(2) \cong \cS^2$. Right-multiplying a rotation matrix $R \in \SO(3)$ by any axial rotation $R_\alpha^{e_3} \in \SO(2)$ leaves the image of the North pole $v = R e_3 \in \cS^2$ unchanged, so the identifiable statistic is precisely the modal axis $v^* = R^* e_3 \in \cS^2$.
\item[-] \textbf{Example~4 (Covariance)}: With $G = \GL(d)$ and $H = N = \Or(d)$, the quotient manifold is $\cM = \GL(d)/\Or(d) \cong \mathrm{SPD}(d)$, the cone of $d \times d$ symmetric positive-definite matrices. Two invertible matrices $A, A' \in \GL(d)$ belong to the same left coset $A\,\Or(d)$ if and only if $A' = AU$ for some $U \in \Or(d)$, which is equivalent to $A'(A')^\top = A A^\top = \Sigma$. Thus, the factor $A^*$ is identifiable up to right-multiplication by $\Or(d)$ (with continuous gradient optimization operating within a single connected component such as $\GL^+(d)$, where the connected stabilizer fiber is $\GL^+(d) \cap \Or(d) = \SO(d)$), and the identifiable parameter on $\cM$ is the covariance matrix $\Sigma^* = A^*(A^*)^\top$.
\end{itemize}

\paragraph{3. Partial overlap ($\{e\} \subsetneq N \subsetneq H$).}
When $G$ and $H$ share a proper, non-trivial subgroup $N = G \cap H$, the identifiable parameter space is the quotient $\cM = G/N$, while the remaining symmetries in $H \setminus N$ act non-trivially on $\cM$ from the outside via conjugation:
\begin{itemize}[leftmargin=*]
\item[-] \textbf{Example~5 (Square)}: With $G = \SO(2)$ and $H = D_4$, the intersection is the discrete fourfold rotation group $N = C_4 = \{R_{k\pi/2} : k = 0,1,2,3\}$. Because rotating a square-symmetric distribution by $90^\circ$ produces the exact same distribution, the orientation angle $\theta \in [0, 2\pi)$ is identifiable only modulo $\pi/2$, on the quotient manifold $\cM = \SO(2)/C_4 \cong \cS^1$. Meanwhile, the reflections in $D_4 \setminus C_4$ do not belong to $G$, yet they act on $\cM$ via conjugation ($\theta \mapsto -\theta$), which will prove decisive in forcing the gradient at $\theta^*$ to vanish.
\end{itemize}

\textbf{Remark (Identifiability versus unbiased stationarity).} Because $q_{gn} = q_g$ for all $n \in N$, the variational loss satisfies $\cL(gn) = \cL(g)$ for every $g \in G$ and $n \in N$. When $N$ is a positive-dimensional Lie subgroup (as in Examples~3 and~4), $\cL$ is identically constant along the coset fibers $gN$, so the directional derivative of $\cL$ at $g^*$ automatically vanishes along all tangent directions to $N$. When $N$ is finite (as in Example~5), $\cL$ is periodic across the $|N|$ equivalent representatives $\{g^* n : n \in N\}$. However, this fiber invariance only accounts for directions \emph{inside} $N$; whether the gradient vanishes along the identifiable directions transverse to $N$ depends on the adjoint action of the full symmetry group $H$ on the Lie algebra of $G$, which we develop in the next two sections.

%% ===========================================================================
\section{Differentiation on the transport Lie group}
\label{sec:lie}
%% ===========================================================================

To characterize when the true transport parameter $g^* \in G$ is a first-order stationary point of the variational loss $\cL(g) = \KL(q_g \| p)$, we require differential calculus on the Lie group $G$. In classical parametric statistics, one typically differentiates a loss function by adding a small coordinate vector $\theta + t u$ in $\R^p$. For transformation parameters such as rotation matrices $R \in \SO(d)$ or scale matrices $A \in \GL(d)$, however, an additive perturbation $A + t U$ either leaves the manifold of valid parameters (breaking orthogonality) or depends on the current scale $A$. Because $G$ is both a group and a smooth manifold, local perturbations around any parameter $g \in G$ can instead be parameterized canonically through the tangent space at the identity---the Lie algebra $\gA = T_e G$---via left multiplication $t \mapsto g \exp(t u)$ (where $\exp(tu) \approx I + tu$ for matrix groups). Writing local perturbations multiplicatively as $g \cdot \exp(tu)$ rather than additively (known geometrically as \emph{left-trivialization}) measures \emph{relative} (or equivariant) parameter changes in a single fixed vector space $\gA$ at the identity. From a statistical perspective, this multiplicative parameterization is natural for transformation models: perturbing $g^*$ to $g^* \exp(tu)$ corresponds, after pulling back by $(g^*)^{-1}$, to perturbing the standardized base distribution $q_0$ by $\exp(tu)$ near the identity $e_G$. We briefly review the necessary tools from Lie theory and Riemannian optimization; see \citet{hall2015lie} and \citet{boumal2023introduction} for comprehensive treatments.

\subsection{Lie group calculus and left-trivialized Riemannian gradients}

\begin{definition}[Lie algebra, exponential map, and Riemannian gradient] \label{def:lie-gradient}
Let $G$ be a Lie group with identity element $e_G$. We recall the following differential-geometric objects:
\begin{itemize}[leftmargin=*]
    \item[-] The \emph{Lie algebra} $\gA \coloneqq T_{e_G}G$ is the tangent vector space to $G$ at the identity, which we equip with a positive-definite inner product $\langle \cdot, \cdot\rangle_\gA$ (defining a left-invariant Riemannian metric on $G$).
    \item[-] The \emph{exponential map} $\exp: \gA \to G$ maps each tangent vector $u \in \gA$ to the time-$1$ point $\gamma_u(1)$ of the unique one-parameter subgroup $\gamma_u(t) = \exp(tu)$ satisfying $\gamma_u(0) = e_G$ and $\gamma_u'(0) = u$. The map $\exp$ is a smooth local diffeomorphism from a neighborhood of $0 \in \gA$ onto a neighborhood of $e_G \in G$, so a neighborhood of any $g \in G$ is parameterized via $u \mapsto g \cdot \exp(u)$.
    \item[-] The (left-trivialized) \emph{Riemannian gradient} $\nabla_G \cL(g) \in \gA$ of a smooth objective $\cL: G \to \R$ at $g \in G$ is the unique element of $\gA$ representing the first-order Taylor expansion of $u \mapsto \cL(g \cdot \exp(u))$ at $u = 0$:
    \begin{equation}\label{eq:riemannian-gradient}
    \cL(g \cdot \exp(u)) = \cL(g) + \langle \nabla_G \cL(g), u \rangle_{\gA} + O(\|u\|^2) \quad \forall u \in \gA \, .
    \end{equation}
\end{itemize}
\end{definition}

\paragraph{Concrete form of the exponential map and Riemannian gradient.}
In statistical applications, $G$ is either an additive vector space or a matrix Lie group, where both the exponential map and the left-trivialized gradient take explicit forms:
\begin{itemize}[leftmargin=*]
\item[-] \textbf{Additive translation group} $G = (\R^d, +)$: The group operation is vector addition, the identity is $e_G = 0$, the Lie algebra is $\gA = \R^d$ equipped with the standard Euclidean inner product, and the exponential map is the identity $\exp(u) = u$. Perturbing $\mu \in G$ via $\mu \cdot \exp(tu)$ is simply $\mu + tu$, and the Riemannian gradient coincides with the ordinary Euclidean gradient:
\[
\nabla_G \cL(\mu) = \nabla_{\text{Euc}} \cL(\mu) \in \R^d\, .
\]
\item[-] \textbf{Matrix Lie groups} $G \subseteq \GL(d)$: When $G$ is a closed subgroup of $\GL(d)$, $\exp(U) = \sum_{k=0}^\infty \frac{1}{k!} U^k$ is the standard matrix exponential, and $g \cdot \exp(tU) = g e^{tU}$ is matrix multiplication. If $\cL$ is extended smoothly to an open neighborhood in $\R^{d \times d}$ with Euclidean matrix of partial derivatives $\nabla_{\text{Euc}} \cL(g) \in \R^{d \times d}$, then for any $U \in \gA$, the chain rule gives
\[
\left.\frac{\ud}{\ud t}\right|_{t=0} \cL(g e^{tU}) = \langle \nabla_{\text{Euc}} \cL(g), \, g U \rangle_{\text{Frob}} = \langle g^\top \nabla_{\text{Euc}} \cL(g), \, U \rangle_{\text{Frob}}\, .
\]
Consequently, under the Frobenius inner product $\langle U, V \rangle_\gA = \mathrm{tr}(U^\top V)$, the Riemannian gradient in $\gA$ is obtained by left-multiplying the Euclidean gradient by $g^\top$ (transporting the tangent vector from $T_g G$ back to $T_I G = \gA$) and projecting orthogonally onto the subspace $\gA \subseteq \R^{d \times d}$:
\[
\nabla_G \cL(g) = \mathrm{proj}_\gA\!\left(g^\top \nabla_{\text{Euc}} \cL(g)\right)\, .
\]
For the general linear group $G = \GL(d)$, the Lie algebra is the full matrix space $\gA = \mathfrak{gl}(d) = \R^{d \times d}$, so $\mathrm{proj}_\gA$ is the identity and $\nabla_G \cL(A) = A^\top \nabla_{\text{Euc}} \cL(A)$. For the special orthogonal group $G = \SO(d)$, the Lie algebra $\gA = \mathfrak{so}(d) = \{U \in \R^{d \times d} : U^\top = -U\}$ consists of skew-symmetric matrices, and orthogonal projection onto $\mathfrak{so}(d)$ is the skew-symmetrization $\mathrm{proj}_{\mathfrak{so}(d)}(M) = \tfrac{1}{2}(M - M^\top)$.
\end{itemize}

\subsection{Adjoint representation and the trivial subspace $\gA^H$}

In Section~\ref{sec:conjugation-detail}, we saw that the symmetry group $H$ acts on the transport group $G$ by conjugation automorphisms $\Phi_h(g) = h \circ g \circ h^{-1}$. While $\Phi_h : G \to G$ is in general a nonlinear map on the manifold $G$, it fixes the identity element ($\Phi_h(e_G) = e_G$ for every $h \in H$). Differentiating $\Phi_h$ at $e_G$ therefore induces a \emph{linear representation} of $H$ on the tangent space $\gA = T_{e_G} G$: the \emph{adjoint representation}.

\begin{definition}[Adjoint representation and trivial subspace] \label{def:adjoint}
The \emph{adjoint representation} of $H$ on $\gA$ is the group homomorphism $\Ad: H \to \GL(\gA)$, $h \mapsto \Ad_h$, defined as the differential of the conjugation map $\Phi_h$ at the identity $e_G$. For any $h \in H$ and $u \in \gA$:
\[
\Ad_h \cdot u \coloneqq (d\Phi_h)_{e_G}(u) = \left.\frac{\ud}{\ud t}\right|_{t=0} \Phi_h(\exp(tu)) \in \gA\, .
\]
The \emph{trivial representation} (or \emph{$H$-fixed subspace}) $\gA^H \subseteq \gA$ is the linear subspace of tangent directions left invariant by the adjoint action of every element of $H$ (called ``trivial'' in representation theory because every $h \in H$ acts on $\gA^H$ as the identity operator, leaving its vectors unchanged):
\[
\gA^H = \{v \in \gA\, :\, \Ad_h \cdot v = v \;\; \forall \, h \in H\}\, .
\]
\end{definition}

\textbf{Compactness of $H$ and $\Ad_H$-invariant inner products.} Because $H$ is assumed to be a compact topological group, it admits a unique normalized Haar probability measure. By averaging any initial inner product $\langle \cdot, \cdot \rangle_0$ on $\gA$ over $H$ via $\langle u, v \rangle_\gA \coloneqq \int_H \langle \Ad_h u, \Ad_h v \rangle_0 \, \ud\mathrm{Haar}(h)$, we are guaranteed an $\Ad_H$-invariant inner product on $\gA$ satisfying
\[
\langle \Ad_h u, \Ad_h v \rangle_\gA = \langle u, v \rangle_\gA \qquad \forall h \in H, \;\; \forall u, v \in \gA\, .
\]
With respect to this inner product, each linear operator $\Ad_h : \gA \to \gA$ is orthogonal ($\Ad_h^\top = \Ad_h^{-1} = \Ad_{h^{-1}}$). Note that for all of our guiding examples, the standard Euclidean inner product on $\RR^d$ and the Frobenius inner product on $\mathfrak{gl}(d)$ or $\mathfrak{so}(d)$ are already $\Ad_H$-invariant without requiring explicit averaging.

\paragraph{Statistical and geometric interpretation of $\Ad_h$ and $\gA^H$.}
Just as $\Phi_h(g)$ represents the transport $g$ viewed in an $h$-transformed coordinate frame, $\Ad_h \cdot u$ represents the infinitesimal parameter perturbation $u \in \gA$ rotated by the symmetry $h$. In our two primary classes of transport groups, $\Ad_h$ takes a particularly transparent form:
\begin{itemize}[leftmargin=*]
\item[-] \emph{Translation groups} ($G = (\RR^d, +)$ with $H$ acting linearly by $h \cdot x = M_h x$): Because $\Phi_h(\mu) = M_h \mu$ is already linear in $\mu$, its derivative at $0$ is simply matrix-vector multiplication:
\[
\Ad_h \cdot u = M_h u \qquad (u \in \RR^d)\, .
\]
\item[-] \emph{Matrix groups} ($G \subseteq \GL(d)$ with $H$ acting linearly by $M_h \in \Or(d)$): Because $\Phi_h(g) = M_h g M_h^{-1}$ is matrix conjugation, differentiating $t \mapsto M_h e^{tU} M_h^{-1}$ at $t = 0$ yields matrix similarity on the Lie algebra:
\[
\Ad_h \cdot U = M_h U M_h^{-1} = M_h U M_h^\top \qquad (U \in \gA \subseteq \R^{d \times d})\, .
\]
\end{itemize}
The trivial subspace $\gA^H = \{u \in \gA : \Ad_h \cdot u = u \;\forall h \in H\}$ is precisely the Lie algebra of the fixed-point subgroup $G^H = \{g \in G : \Phi_h(g) = g \;\forall h \in H\}$ introduced in Section~\ref{sec:conjugation-detail}. Statistically, a tangent vector $u$ belongs to $\gA^H$ if and only if the one-parameter flow $\exp(tu)$ commutes with every symmetry $h \in H$, which means that pushing an $H$-symmetric distribution $q_0$ forward along $\exp(tu)$ \emph{preserves its $H$-symmetry}. Conversely, any direction $u$ orthogonal to $\gA^H$ breaks the $H$-symmetry of $q_0$.

\paragraph{Guiding examples: Adjoint action and trivial subspace $\gA^H$.}
We now determine $\gA$, the linear operator $\Ad_h$, and the trivial subspace $\gA^H$ for each of the six guiding examples:

\begin{enumerate}[leftmargin=*]
\item \textbf{Reflection} ($G = (\R, +)$, $H = \Z/2\Z = \{1, -1\}$):
The Lie algebra is $\gA = \R$, and the non-trivial reflection $\sigma = -1$ acts on tangent shifts $u \in \R$ by sign inversion, $\Ad_\sigma u = -u$. The only shift direction satisfying $-u = u$ is $u = 0$, so the trivial subspace vanishes: $\gA^H = \{0\}$. Every nonzero translation breaks reflection symmetry around the origin.

\item \textbf{Euclidean} ($G = (\R^d, +)$, $H = \Or(d)$):
The Lie algebra is $\gA = \R^d$, and each orthogonal matrix $U \in \Or(d)$ acts on tangent vectors $u \in \R^d$ by $\Ad_U u = Uu$. Because no nonzero vector in $\R^d$ is invariant under all orthogonal transformations, the trivial subspace is again trivial: $\gA^H = \{0\}$.

\item \textbf{Sphere} ($G = \SO(3)$, $H = \SO(2)$ fixing $e_3$):
The Lie algebra $\gA = \mathfrak{so}(3)$ of $3 \times 3$ skew-symmetric matrices can be identified with $\R^3$ via the hat isomorphism $\omega \in \R^3 \leftrightarrow [\omega]_\times \in \mathfrak{so}(3)$, where $[\omega]_\times x = \omega \times x$. Under this identification, matrix conjugation $\Ad_h([\omega]_\times) = R_\alpha^{e_3} [\omega]_\times (R_\alpha^{e_3})^\top = [R_\alpha^{e_3} \omega]_\times$ acts on the instantaneous rotation axis $\omega \in \R^3$ simply by rotating $\omega$ around $e_3$ by angle $\alpha$. The vectors invariant under all rotations around $e_3$ are those parallel to $e_3$ itself, yielding a one-dimensional trivial subspace:
\[
\gA^H = \mathrm{span}([e_3]_\times) \cong \mathrm{span}(e_3)\, .
\]
Note that this subspace coincides with the Lie algebra $\mathfrak{n} = \mathfrak{so}(2)$ of the stabilizer $N = \SO(2)$ of $q_0$.

\item \textbf{Covariance} ($G = \GL(d)$, $H = \Or(d)$):
The Lie algebra is the space of all $d \times d$ matrices $\gA = \mathfrak{gl}(d) = \R^{d \times d}$, on which $U \in \Or(d)$ acts by orthogonal similarity $\Ad_U(M) = U M U^\top$. A matrix $M \in \mathfrak{gl}(d)$ belongs to $\gA^H$ if and only if $U M U^\top = M$ (i.e., $U M = M U$) for all $U \in \Or(d)$. By Schur's lemma, the only matrices commuting with all orthogonal transformations are scalar multiples of the identity matrix:
\[
\gA^H = \mathrm{span}(I_d) = \{\lambda I_d : \lambda \in \R\}\, .
\]
Statistically, $\mathrm{span}(I_d)$ generates the one-parameter subgroup $\exp(t I_d) = e^t I_d$ of isotropic scalings, which dilate or contract the overall variance of an $\Or(d)$-symmetric distribution without breaking its spherical symmetry.

\item \textbf{Square} ($G = \SO(2)$, $H = D_4$):
The Lie algebra of planar rotations is one-dimensional, $\gA = \mathfrak{so}(2) \cong \R$, parameterized by the angular velocity $\dot{\theta} \in \R$. Differentiating $\Phi_{r_k}(R_\theta) = R_\theta$ and $\Phi_{s_k}(R_\theta) = R_{-\theta}$ at $\theta = 0$ shows that rotations $r_k \in C_4$ act trivially ($\Ad_{r_k} \dot{\theta} = \dot{\theta}$), whereas reflections $s_k \in D_4 \setminus C_4$ reverse the angular velocity ($\Ad_{s_k} \dot{\theta} = -\dot{\theta}$). Because an invariant direction must be fixed by \emph{all} elements of $D_4$ including reflections, we obtain $\gA^H = \{0\}$.

\item \textbf{Permutation} ($G = (\R^n, +)$, $H = \Sn$):
The Lie algebra is $\gA = \R^n$, and each permutation $\sigma \in \Sn$ acts on tangent shift vectors $u \in \R^n$ by permuting their coordinates, $\Ad_\sigma u = P_\sigma u$. Under the standard representation of $\Sn$, $\R^n$ decomposes orthogonally into the one-dimensional trivial representation of constant vectors and the $(n-1)$-dimensional standard irreducible representation of zero-sum vectors:
\[
\gA^H = \{u \in \R^n : P_\sigma u = u \;\; \forall \sigma \in \Sn\} = \mathrm{span}(\ones)\, , \qquad (\gA^H)^\perp = \Big\{u \in \R^n : \sum_{i=1}^n u_i = 0\Big\}\, .
\]
Statistically, perturbations along $\mathrm{span}(\ones)$ shift all coordinates by a common scalar additive constant, preserving exchangeability, whereas perturbations in $(\gA^H)^\perp$ alter the relative offsets between coordinates and break exchangeability.
\end{enumerate}

As we prove in Section~\ref{sec:main-result}, the Riemannian gradient $\nabla_G \cL(g^*)$ is forced to lie inside $\gA^H$. Whenever $\gA^H = \{0\}$ (Examples~1, 2, and~5), first-order stationarity $\nabla_G \cL(g^*) = 0$ follows immediately from symmetry alone; when $\gA^H \neq \{0\}$ (Examples~3, 4, and~6), the interplay between $\gA^H$ and the identifiability kernel $\mathfrak{n}$ determines whether any nonzero gradient component can survive.

%% ===========================================================================
\section{Main results: symmetry constraints, Lie algebra decomposition, and moment matching} \label{sec:main-result}
%% ===========================================================================

We now establish our main theoretical results characterizing the first-order optimization landscape of the variational objective $\cL(g) = \KL(q_g \| p)$ at the true transport parameter $g^*$. Geometrically, our results rest on a fundamental principle of invariant analysis: \emph{the gradient of a group-invariant functional at a symmetric point must itself be invariant under the group action}. To picture this concretely, consider a rotationally invariant (radial) function on $\R^d$: computing its gradient at a rotated point on the same circular level set yields an equally rotated gradient vector, and at the center of symmetry itself (the origin), the gradient vector must remain unchanged under every rotation---which forces it to be the zero vector. When translated to the Lie group setting via the pullback loss around $g^*$, this principle implies that the Riemannian gradient $\nabla_G \cL(g^*)$ is fixed by the adjoint action $\Ad_h$ for every $h \in H$, and therefore must reside in the trivial representation $\gA^H$.

\subsection{Warm-up: Gradient constraint for translation families}

To build concrete statistical intuition before treating general Lie groups, we first state and prove the gradient constraint for additive location families $G = (\R^d, +)$, where the Lie algebra is $\R^d$ and the Riemannian gradient reduces to the standard Euclidean gradient. This proposition generalizes the univariate reflection argument of \citet{margossian2024variational} (where $d = 1$ and $H = \{\pm 1\}$) to arbitrary linear isometry groups $H \subseteq \Or(d)$ in $d$ dimensions, and illustrates how the expected score of the standardized target interacts with $H$-invariance.

\begin{proposition}[Gradient constraint --- translation case] \label{prop:gradient-translation}
In the setting of Definition~\ref{def:vi-symmetry}, suppose that $G = (\R^d, +)$ acts on $\cX = \R^d$ by translation $g \cdot x = x + \mu_g$, and that $H$ acts linearly by orthogonal transformations $h \cdot x = M_h x$. Then the Euclidean gradient of $\cL(\mu) = \KL(q_\mu \| p)$ at the true location parameter $\mu^*$ satisfies
\[
\nabla \cL(\mu^*) = M_h \nabla \cL(\mu^*)\,, \quad \forall h \in H\, .
\]
Equivalently, $\nabla \cL(\mu^*)$ lies in the $H$-fixed subspace $(\R^d)^H = \{v \in \R^d : M_h v = v \;\forall h \in H\}$.
\end{proposition}

\begin{proof}
Applying the translation change of variables $x = \mu + \xi$ (with unit Jacobian determinant), the variational objective at any location $\mu \in \R^d$ can be written as an expectation under the fixed base distribution $\xi \sim q_0$:
\begin{align*}
\cL(\mu) &= \int_{\R^d} q_0(\xi) \big[\log q_0(\xi) - \log p(\mu + \xi)\big] \, \ud \xi\, .
\end{align*}
Differentiating with respect to $\mu$ under the integral sign and evaluating at $\mu = \mu^*$, where $p(\mu^* + \xi) = \tilde p(\xi)$ by definition of the standardized pullback $\tilde p$, yields the expected score of $\tilde p$ under $q_0$:
\begin{align*}
\nabla_G\cL(\mu^*) = \nabla \cL(\mu^*) &= -\int_{\R^d} q_0(\xi) \nabla_\mu\big[\log p(\mu^* + \xi)\big] \, \ud \xi = -\int_{\R^d} q_0(\xi) \nabla \log \tilde p(\xi)\, \ud \xi \, .
\end{align*}
Because $\tilde p$ is $H$-symmetric ($\tilde p(M_h \xi') = \tilde p(\xi')$), the chain rule implies that its score function is $H$-equivariant: $M_h^\top \nabla \log \tilde p(M_h \xi') = \nabla \log \tilde p(\xi')$, or equivalently $\nabla \log \tilde p(M_h \xi') = M_h \nabla \log \tilde p(\xi')$ since $M_h \in \Or(d)$ is orthogonal. Performing the orthogonal change of variables $\xi = M_h \xi'$ (with $|\det M_h| = 1$) and using the $H$-invariance $q_0(M_h \xi') = q_0(\xi')$, we obtain for every $h \in H$:
\begin{align*}
\nabla \cL(\mu^*) &= -\int_{\R^d} q_0(M_h \xi') \, \nabla \log \tilde p(M_h \xi') \, \ud \xi' = -M_h \int_{\R^d} q_0(\xi') \, \nabla \log \tilde p(\xi') \, \ud \xi' = M_h \nabla \cL(\mu^*)\, .
\end{align*}
Thus $v^* = \nabla \cL(\mu^*)$ satisfies $M_h v^* = v^*$ for all $h \in H$, i.e., $\nabla \cL(\mu^*) \in (\R^d)^H$. In Examples~1 (Reflection) and~2 (Euclidean), $(\R^d)^H = \{0\}$, which immediately forces $\nabla \cL(\mu^*) = 0$; in Example~6 (Permutation), it constrains $\nabla \cL(\mu^*) \in \mathrm{span}(\ones)$.
\end{proof}

\subsection{Gradient constraint: the general Lie group case}

We now establish our first main theorem for an arbitrary transport Lie group $G$ normalized by a compact symmetry group $H$. In the translation setting of Proposition~\ref{prop:gradient-translation}, the matrix $M_h$ acted directly on the parameter gradient because the conjugation of a translation $\mu$ by $M_h$ is simply $\Phi_h(\mu) = M_h \mu$. For a general Lie group $G$, the natural counterpart of $M_h$ on the tangent space $\gA = T_e G$ is the adjoint operator $\Ad_h = (d\Phi_h)_{e_G}$. Replacing Euclidean invariance with the conjugation invariance of the relative loss established in Proposition~\ref{prop:conjugation-distribution}(iii) yields the general gradient constraint.

\begin{theorem}[Gradient constraint from symmetry] \label{thm:gradient-constraint}
Under the assumptions of Definition~\ref{def:vi-symmetry}, let $\gA$ be equipped with an $\Ad_H$-invariant inner product $\langle \cdot, \cdot \rangle_\gA$. Then the Riemannian gradient of the variational objective $\cL(g) = \KL(q_g \| p)$ at the true transport parameter $g^* \in G$ satisfies
\[
\nabla_G\cL(g^*) = \Ad_h \cdot \nabla_G\cL(g^*)\,, \quad \forall h \in H\, .
\]
Equivalently, $\nabla_G\cL(g^*)$ is constrained to lie in the trivial representation of $H$ in $\gA$:
\[
\nabla_G\cL(g^*) \in \gA^H\, .
\]
\end{theorem}

\begin{proof}
For any tangent vector $u \in \gA$, consider the perturbed parameter $g^* \exp(u) \in G$. Because the Kullback--Leibler divergence is invariant under simultaneous application of the diffeomorphism $(g^*)^{-1}$ to both measures, we can pull the variational objective back to the standardized reference space around the identity $e_G$:
\[
\cL(g^* \exp(u)) = \KL\big(q_{g^* \exp(u)} \,\big\|\, p\big) = \KL\big(q_{\exp(u)} \,\big\|\, \tilde p\big) \eqqcolon \cL_0(\exp(u))\, .
\]
By Definition~\ref{def:lie-gradient}, the left-trivialized Riemannian gradient of $\cL$ at $g^*$ coincides with the Riemannian gradient of the standardized loss $\cL_0$ at $e_G$, and hence with the Euclidean gradient at $0 \in \gA$ of the pullback objective $J : \gA \to \R$ defined by $J(u) \coloneqq \cL_0(\exp(u)) = \cL(g^* \exp(u))$.

We first show that the function $J : \gA \to \R$ is strictly invariant under the adjoint action of $H$, i.e., $J(\Ad_h \cdot u) = J(u)$ for all $h \in H$ and $u \in \gA$:
\begin{align*}
    J(\Ad_h \cdot u) &= \cL_0\big(\exp(\Ad_h\cdot u)\big) = \cL_0\big(\Phi_h(\exp(u))\big) \quad \text{($\dagger$)}\\
    &= \KL\big(q_{\Phi_h(\exp(u))} \,\big\|\, \tilde p\big) = \KL\big(q_{h\exp(u)} \,\big\|\, \tilde p\big) \quad \text{(*)}\\
    &= \KL\big(q_{h\exp(u)} \,\big\|\, h_\sharp \tilde p\big) \quad \text{(**)}\\
    &= \KL\big(q_{\exp(u)} \,\big\|\, \tilde p\big) = J(u)\, .
\end{align*}
Each step is justified by the algebraic and measure-theoretic structure established earlier:
\begin{itemize}[leftmargin=*]
\item[-] \textbf{Step ($\dagger$) (Naturality of the exponential map):} $\exp(\Ad_h \cdot u) = \Phi_h(\exp(u))$. Since $\Phi_h : G \to G$ is a Lie group automorphism (Proposition~\ref{prop:conjugation-homomorphism}), the curve $t \mapsto \Phi_h(\exp(tu))$ is a one-parameter subgroup of $G$ whose initial velocity at $t = 0$ is $(d\Phi_h)_{e_G}(u) = \Ad_h \cdot u$ by Definition~\ref{def:adjoint}; by uniqueness of one-parameter subgroups, $\Phi_h(\exp(tu)) = \exp(t \Ad_h \cdot u)$.
\item[-] \textbf{Step (*) (Pushforward equivariance of $q_0$):} $q_{\Phi_h(\exp(u))} = q_{h\exp(u)} = h_\sharp q_{\exp(u)}$, which follows directly from Proposition~\ref{prop:conjugation-distribution}(i) using the $H$-invariance of $q_0$.
\item[-] \textbf{Step (**) ($H$-invariance of the pullback $\tilde p$):} Because $p$ has $H$-conjugate symmetry around $g^*$, its pullback $\tilde p$ satisfies $h_\sharp \tilde p = \tilde p$, and the final equality follows by applying the diffeomorphism $h^{-1}$ to both arguments of the KL divergence.
\end{itemize}

Differentiating the invariance identity $J(\Ad_h \cdot u) = J(u)$ at $u = 0$ along an arbitrary direction $u \in \gA$ via the first-order expansion~\eqref{eq:riemannian-gradient} gives
\[
\langle \nabla_G \cL(g^*), u \rangle_\gA = \langle \nabla J(0), u \rangle_\gA = \langle \nabla J(0), \Ad_h \cdot u \rangle_\gA = \langle \Ad_h^\top \nabla J(0), u \rangle_\gA \qquad \forall u \in \gA\, .
\]
Because the inner product $\langle \cdot, \cdot \rangle_\gA$ is $\Ad_H$-invariant, $\Ad_h$ is an orthogonal operator on $\gA$, so $\Ad_h^\top = \Ad_h^{-1} = \Ad_{h^{-1}}$ (as is also immediately verified in coordinates: for translation groups $\Ad_h^\top = M_h^\top = M_h^{-1}$, and for matrix groups under the Frobenius inner product $\langle \Ad_h U, V \rangle = \mathrm{tr}((M_h U M_h^\top)^\top V) = \mathrm{tr}(U^\top M_h^\top V M_h) = \langle U, \Ad_{h^{-1}} V \rangle$). Since $u \in \gA$ is arbitrary, we conclude that $\nabla_G \cL(g^*) = \Ad_{h^{-1}} \cdot \nabla_G \cL(g^*)$ for all $h \in H$. Replacing $h$ by $h^{-1}$ yields $\nabla_G \cL(g^*) \in \gA^H$.
\end{proof}

\subsection{Gradient localization and orthogonal Lie algebra decomposition}

Theorem~\ref{thm:gradient-constraint} establishes that the Riemannian gradient $\nabla_G \cL(g^*)$ is confined to the trivial representation $\gA^H$. At the same time, our identifiability analysis in Section~\ref{sec:identifiability} showed that whenever the intersection subgroup $N = G \cap H$ has positive dimension, the variational objective $\cL(g^* \exp(tu)) = \cL(g^*)$ is identically constant along the Lie algebra $\mathfrak{n} = T_e N$ of $N$. Combining these two independent structural constraints restricts the gradient to the orthogonal complement of $\mathfrak{n}$ inside $\gA^H$, and induces a canonical three-way orthogonal decomposition of the entire transport Lie algebra $\gA$. From a statistical estimation standpoint, this decomposes all local parameter perturbations around $g^*$ into three intuitive categories: (1) \emph{unidentifiable directions} $\mathfrak{n}$ that do not change the distribution $q_{g^*}$ at all; (2) \emph{symmetry-preserving directions} $\gA^H_{\text{bias}}$ (such as global scale) where mismatch between $q_0$ and $\tilde p$ can induce estimation bias; and (3) \emph{symmetry-breaking directions} $\gA^\perp$ (such as location, orientation, or correlation shape) where symmetry forces the gradient at $g^*$ to vanish.

\begin{proposition}[Gradient localization and Lie algebra decomposition] \label{prop:decomposition}
Under the setting of Definition~\ref{def:vi-symmetry} with an $\Ad_H$-invariant inner product on $\gA$, the transport Lie algebra $\gA$ admits the orthogonal direct sum decomposition
\[
\gA = \underbrace{\mathfrak{n}}_{\text{identifiability kernel (flat)}} \;\oplus\; \underbrace{\gA^H_{\text{bias}}}_{\text{bias directions}} \;\oplus\; \underbrace{\gA^\perp}_{\text{symmetry-breaking (forced to }0\text{)}}\, ,
\]
where the three mutually orthogonal subspaces are defined as follows:
\begin{itemize}[leftmargin=*]
\item[-] \textbf{Identifiability kernel} $\mathfrak{n} \coloneqq \{ u \in \gA \mid \exp(tu) \in N = G \cap H \;\; \forall t \in \R\}$:
the Lie algebra of the stabilizer intersection $N = G \cap H$. Along any direction $u \in \mathfrak{n}$, the transported distribution $q_{g^* \exp(tu)} = q_{g^*}$ is unchanged, so $\cL(g^* \exp(tu)) = \cL(g^*)$ is locally constant and $\mathrm{proj}_{\mathfrak{n}}\big(\nabla_G \cL(g^*)\big) = 0$.
\item[-] \textbf{Bias subspace} $\gA_{\text{bias}}^H \coloneqq \gA^H \cap \mathfrak{n}^\perp$:
the subspace of $H$-invariant tangent directions orthogonal to the identifiability kernel $\mathfrak{n}$. Perturbations along $\gA_{\text{bias}}^H$ deform the distribution $q_{g^*}$ while preserving its $H$-conjugate symmetry, and represent the only directions along which $\nabla_G \cL(g^*)$ can be nonzero.
\item[-] \textbf{Symmetry-breaking subspace} $\gA^\perp \coloneqq (\gA^H)^\perp \cap \mathfrak{n}^\perp$:
the subspace of identifiable tangent directions that transform non-trivially under $\Ad_H$. By Theorem~\ref{thm:gradient-constraint}, the gradient vanishes identically on this subspace: $\mathrm{proj}_{\gA^\perp}\big(\nabla_G \cL(g^*)\big) = 0$.
\end{itemize}
Consequently, the Riemannian gradient at the true transport parameter $g^*$ satisfies
\[
\nabla_G \cL(g^*) \in \gA^H_{\text{bias}}\, .
\]
\end{proposition}

\begin{proof}
By Theorem~\ref{thm:gradient-constraint}, $\nabla_G\cL(g^*) \in \gA^H$. Furthermore, for any $u \in \mathfrak{n}$ and any $t \in \R$, we have $\exp(tu) \in N = G \cap H \subseteq H$. Because $q_0$ is $H$-invariant, $q_{\exp(tu)} = (\exp(tu))_\sharp q_0 = q_0$, which implies $q_{g^*\exp(tu)} = q_{g^*}$ and hence $\cL(g^*\exp(tu)) = \cL(g^*)$ for all $t \in \R$. Differentiating at $t = 0$ gives $\langle \nabla_G \cL(g^*), u \rangle_\gA = 0$ for all $u \in \mathfrak{n}$, so $\nabla_G \cL(g^*) \in \mathfrak{n}^\perp$. Combining both inclusions yields $\nabla_G\cL(g^*) \in \gA^H \cap \mathfrak{n}^\perp = \gA^H_{\text{bias}}$.

To prove that $\gA = \mathfrak{n} \oplus \gA^H_{\text{bias}} \oplus \gA^\perp$ is an orthogonal direct sum spanning $\gA$, note first that $\gA^H_{\text{bias}} \subseteq \mathfrak{n}^\perp$ and $\gA^\perp \subseteq \mathfrak{n}^\perp$ are both orthogonal to $\mathfrak{n}$ by definition, and $\gA^H_{\text{bias}} \subseteq \gA^H$ is orthogonal to $\gA^\perp \subseteq (\gA^H)^\perp$. Next, observe that $N = G \cap H$ is normalized by $H$: for any $n \in N$ and $h \in H$, $\Phi_h(n) = hnh^{-1}$ lies in $G$ (since $G$ is normalized by $H$) and in $H$ (since $H$ is a group), so $\Phi_h(N) = N$. Differentiating at $e_G$ implies that $\mathfrak{n}$ is an $\Ad_H$-invariant subspace of $\gA$ ($\Ad_h(\mathfrak{n}) = \mathfrak{n}$ for all $h \in H$; see Appendix~\ref{app:centrality} for details). Because $\mathfrak{n}$ is $\Ad_H$-invariant and the inner product is $\Ad_H$-invariant, the orthogonal projection onto $\gA^H$ maps $\mathfrak{n}$ into itself, so $\mathfrak{n}$ decomposes orthogonally as
\[
\mathfrak{n} = (\mathfrak{n} \cap \gA^H) \;\oplus\; (\mathfrak{n} \cap (\gA^H)^\perp)\, .
\]
Taking orthogonal complements inside $\gA^H$ and $(\gA^H)^\perp$ respectively gives $\gA^H = (\mathfrak{n} \cap \gA^H) \oplus \gA^H_{\text{bias}}$ and $(\gA^H)^\perp = (\mathfrak{n} \cap (\gA^H)^\perp) \oplus \gA^\perp$. Summing these decompositions over $\gA = \gA^H \oplus (\gA^H)^\perp$ yields $\gA = \mathfrak{n} \oplus \gA^H_{\text{bias}} \oplus \gA^\perp$.
\end{proof}

\paragraph{Guiding examples: Three-way decomposition and statistical interpretation.}
We now evaluate the orthogonal decomposition $\gA = \mathfrak{n} \oplus \gA^H_{\text{bias}} \oplus \gA^\perp$ for each of our six guiding examples, illustrating how Proposition~\ref{prop:decomposition} separates overparametrization, estimation bias, and unbiased recovery:

\begin{enumerate}[leftmargin=*]
    \item \textbf{Reflection} ($G = (\R, +)$, $H = \Z/2\Z$):
    The Lie algebra is $\gA = \R$. Because $N = \{0\}$ is trivial and $\Ad_\sigma u = -u$ leaves only $0$ fixed, both the identifiability kernel and the bias subspace vanish:
    \[
    \mathfrak{n} = \{0\}\, , \qquad \gA^H_{\text{bias}} = \{0\}\, , \qquad \gA^\perp = \R\, .
    \]
    Every nonzero translation breaks reflection symmetry, so the entire parameter space lies in $\gA^\perp$ and $\nabla_G \cL(\mu^*) = 0$ is guaranteed for any symmetric target $p$.

    \item \textbf{Euclidean} ($G = (\R^d, +)$, $H = \Or(d)$):
    The Lie algebra is $\gA = \R^d$, with $N = \{0\}$ and $\gA^H = \{0\}$. Thus
    \[
    \mathfrak{n} = \{0\}\, , \qquad \gA^H_{\text{bias}} = \{0\}\, , \qquad \gA^\perp = \R^d\, .
    \]
    All $d$ translation coordinates are symmetry-breaking directions, forcing $\nabla_G \cL(\mu^*) = 0$.

    \item \textbf{Sphere} ($G = \SO(3)$, $H = \SO(2)$ fixing $e_3$):
    Under the isomorphism $\gA = \mathfrak{so}(3) \cong \R^3$, the stabilizer $N = \SO(2)$ has Lie algebra $\mathfrak{n} = \mathrm{span}(e_3)$, consisting of infinitesimal rotations around the symmetry axis $e_3$. As computed in Section~\ref{sec:lie}, the trivial representation is also $\gA^H = \mathrm{span}(e_3)$. Because $\gA^H$ is entirely contained within the identifiability kernel $\mathfrak{n}$, the bias subspace vanishes:
    \[
    \mathfrak{n} = \mathrm{span}(e_3)\, , \qquad \gA^H_{\text{bias}} = \gA^H \cap \mathfrak{n}^\perp = \{0\}\, , \qquad \gA^\perp = \mathrm{span}(e_1, e_2)\, .
    \]
    Geometrically, $\gA^\perp = \mathrm{span}(e_1, e_2)$ consists of infinitesimal tilts of the modal axis $e_3$ toward the equatorial plane (spanning the tangent space $T_{e_3}\cS^2$ of the quotient manifold $\cM = \SO(3)/\SO(2) \cong \cS^2$). Symmetry forces the gradient along all equatorial tilts to vanish, while the axial spin direction $\mathfrak{n}$ is flat by identifiability, guaranteeing $\nabla_G \cL(R^*) = 0$.

    \item \textbf{Covariance} ($G = \GL(d)$, $H = \Or(d)$):
    The Lie algebra $\gA = \mathfrak{gl}(d) = \R^{d \times d}$ equipped with the Frobenius inner product decomposes into skew-symmetric matrices and symmetric matrices. The stabilizer $N = \Or(d)$ has Lie algebra $\mathfrak{n} = \mathfrak{so}(d) = \{U \in \R^{d \times d} : U^\top = -U\}$, while the trivial subspace is the one-dimensional space of isotropic scalings $\gA^H = \mathrm{span}(I_d)$. Since $I_d$ is symmetric and every matrix in $\mathfrak{so}(d)$ is skew-symmetric (and traceless), $\gA^H$ is completely orthogonal to $\mathfrak{n}$ ($\mathrm{tr}(I_d^\top U) = \mathrm{tr}(U) = 0$ for all $U \in \mathfrak{so}(d)$). We therefore obtain the classical Cartan-type matrix decomposition:
    \[
    \mathfrak{n} = \mathfrak{so}(d)\, , \qquad \gA^H_{\text{bias}} = \mathrm{span}(I_d)\, , \qquad \gA^\perp = \mathrm{Sym}_0(d) \coloneqq \{S \in \R^{d \times d} : S = S^\top, \, \mathrm{tr}(S) = 0\}\, .
    \]
    This decomposition has a striking statistical interpretation. The identifiability kernel $\mathfrak{n} = \mathfrak{so}(d)$ generates orthogonal right-rotations of $A^*$ that leave the covariance matrix $\Sigma^* = A^*(A^*)^\top$ invariant. The symmetry-breaking subspace $\gA^\perp = \mathrm{Sym}_0(d)$ of traceless symmetric matrices parameterizes \emph{anisotropic shape and correlation deformations} (changing the ratios of eigenvalues and orientation of principal axes while preserving volume to first order); along all $\tfrac{1}{2}d(d+1) - 1$ such directions, the gradient at $A^*$ is forced to vanish by elliptical symmetry! However, along the single remaining scalar direction $\gA^H_{\text{bias}} = \mathrm{span}(I_d)$, the gradient satisfies $\nabla_G \cL(A^*) = \lambda I_d$ for some scalar $\lambda \in \R$, corresponding to a \emph{multiplicative global scale bias} between the dispersion of $q_0$ and that of $\tilde p$.

    \item \textbf{Square} ($G = \SO(2)$, $H = D_4$):
    The Lie algebra is $\gA = \mathfrak{so}(2) \cong \R$. Because the overlap subgroup $N = C_4$ is discrete, its Lie algebra is zero-dimensional ($\mathfrak{n} = \{0\}$). At the same time, the reflections in $D_4$ force $\gA^H = \{0\}$. Consequently,
    \[
    \mathfrak{n} = \{0\}\, , \qquad \gA^H_{\text{bias}} = \{0\}\, , \qquad \gA^\perp = \R\, .
    \]
    Every infinitesimal planar rotation breaks the reflection symmetry of the square, so $\nabla_G \cL(R_{\theta^*}) = 0$ is guaranteed by $D_4$-invariance.

    \item \textbf{Permutation} ($G = (\R^n, +)$, $H = \Sn$):
    The Lie algebra is $\gA = \R^n$ with trivial intersection $N = \{0\}$, so $\mathfrak{n} = \{0\}$. The trivial representation of $\Sn$ is the one-dimensional diagonal subspace $\gA^H = \mathrm{span}(\ones)$, which therefore coincides entirely with the bias subspace:
    \[
    \mathfrak{n} = \{0\}\, , \qquad \gA^H_{\text{bias}} = \mathrm{span}(\ones)\, , \qquad \gA^\perp = \Big\{u \in \R^n : \sum_{i=1}^n u_i = 0\Big\}\, .
    \]
    Statistically, the $(n-1)$-dimensional symmetry-breaking subspace $\gA^\perp$ of zero-sum vectors governs all \emph{relative contrasts} $\mu_i - \mu_j$ between coordinates: exchangeability forces the gradient at $\mu^*$ to vanish along every contrast vector $u \in \gA^\perp$, ensuring that relative coordinate differences are estimated without bias. By contrast, along the common-location direction $\gA^H_{\text{bias}} = \mathrm{span}(\ones)$, the gradient takes the form $\nabla_G \cL(\mu^*) = \lambda \ones$ for some $\lambda \in \R$, representing a uniform \emph{additive location bias} across all coordinates if $q_0$ and $\tilde p$ are not centered at the same scalar mean/median level.
\end{enumerate}

Proposition~\ref{prop:decomposition} unifies the phenomena observed in \citet{margossian2024variational} and extends them across general Lie groups. Whenever $\gA^H_{\text{bias}} = \{0\}$---either because $\gA^H = \{0\}$ (Reflection, Euclidean, Square) or because $\gA^H = \mathfrak{n}$ is absorbed into the stabilizer of $q_0$ (Sphere)---symmetry alone guarantees $\nabla_G \cL(g^*) = 0$. Conversely, when $\gA^H_{\text{bias}} \neq \{0\}$ (Covariance and Permutation), symmetry eliminates all gradient components in $\gA^\perp$ while isolating any potential misspecification bias to the low-dimensional subspace $\gA^H_{\text{bias}}$. Table~\ref{tab:summary} summarizes the algebraic invariants and decomposition across all six guiding examples.

\begin{table}[ht!]
\centering
\caption{Algebraic invariants and Lie algebra decomposition $\gA = \mathfrak{n} \oplus \gA^H_{\text{bias}} \oplus \gA^\perp$ for the six guiding examples. The final column indicates whether $\nabla_G\cL(g^*) = 0$ is guaranteed by symmetry alone ($\gA^H_{\text{bias}} = \{0\}$).}
\label{tab:summary}
{\small
\begin{tabular}{l|c|c|c|c|c|c|c|c}
\textbf{Example} & $\cX$ & $G$ & $H$ & $N = G \cap H$ & Regime & $\gA^H$ & $\gA^H_{\text{bias}}$ & $\nabla = 0$? \\ \hline
1. Reflection & $\R$ & $(\R, +)$ & $\Z/2\Z$ & $\{e\}$ & Separable & $\{0\}$ & $\{0\}$ & Yes \\
2. Euclidean & $\R^d$ & $(\R^d, +)$ & $\Or(d)$ & $\{e\}$ & Separable & $\{0\}$ & $\{0\}$ & Yes \\
3. Sphere & $\cS^2$ & $\SO(3)$ & $\SO(2)$ & $\SO(2)$ & Internal & $\text{span}(e_3)$ & $\{0\}$ & Yes \\
4. Covariance & $\R^d$ & $\GL(d)$ & $\Or(d)$ & $\Or(d)$ & Internal & $\text{span}(I_d)$ & $\text{span}(I_d)$ & No \\
5. Square & $\R^2$ & $\SO(2)$ & $D_4$ & $C_4$ & Partial & $\{0\}$ & $\{0\}$ & Yes \\
6. Permutation & $\R^n$ & $(\R^n, +)$ & $\Sn$ & $\{e\}$ & Separable & $\text{span}(\ones)$ & $\text{span}(\ones)$ & No \\
\end{tabular}
}
\end{table}

\subsection{Bias from model mismatch and the Stein moment-matching condition}

When the bias subspace $\gA^H_{\text{bias}} = \gA^H \cap \mathfrak{n}^\perp$ is non-trivial, Proposition~\ref{prop:decomposition} shows that $\nabla_G \cL(g^*)$ lies in $\gA^H_{\text{bias}}$ but need not vanish. Statistically, a nonzero gradient component along $\gA^H_{\text{bias}}$ arises because perturbations $\exp(tu)$ with $u \in \gA^H_{\text{bias}}$ transform the base distribution $q_0$ into a distinct $H$-symmetric distribution $q_{\exp(tu)}$. If the standardized target $\tilde p$ and the base $q_0$ differ in their scale or position along this symmetry-preserving subgroup $\exp(\gA^H_{\text{bias}})$, the KL divergence pulls the variational parameter away from $g^*$ along $\gA^H_{\text{bias}}$ to compensate for the mismatch between $q_0$ and $\tilde p$.

\paragraph{Illustrative examples of estimation bias under model mismatch.}
Before deriving the general analytical criterion for when the bias vanishes, we illustrate how mismatch along $\gA^H_{\text{bias}}$ manifests in three concrete settings:
\begin{itemize}[leftmargin=*]
\item[-] \textbf{Multiplicative scale bias in covariance estimation (Example~4)}: Suppose the base distribution is $q_0 = \cN(0, I_d)$, while the standardized pullback of the target is $\tilde p = \cN(0, \sigma^2 I_d)$ with $\sigma \neq 1$ (so that $p = \cN(0, \sigma^2 A^*(A^*)^\top)$). Both $q_0$ and $\tilde p$ are $\Or(d)$-symmetric, yet their radial scales differ by $\sigma$. Consequently, $\cL(A)$ is minimized at $\hat A = \sigma A^*$ rather than $A^*$, and the gradient at $A^*$ has a nonzero component $\nabla_G \cL(A^*) = (\sigma^{-2} - 1) I_d \in \mathrm{span}(I_d) = \gA^H_{\text{bias}}$ (so that the negative gradient $-\nabla_G \cL(A^*) = (1 - \sigma^{-2})I_d$ drives $A$ toward $\sigma A^*$).

\item[-] \textbf{Common-location additive bias under exchangeability (Example~6)}: Suppose the base distribution is $q_0(x) \propto \exp(-\kappa\|x\|_1)$ (centered at the origin), while the standardized pullback is $\tilde p(x) \propto \exp(-\kappa\|x - \beta\ones\|_1)$ with $\beta \neq 0$. Because shifting all coordinates by $\beta\ones$ commutes with all coordinate permutations in $\Sn$, $\tilde p$ remains fully $\Sn$-symmetric! However, $q_0$ and $\tilde p$ have different common coordinate medians ($0$ versus $\beta$), so $\cL(\mu)$ is minimized at $\hat \mu = \mu^* + \beta\ones$ rather than $\mu^*$, producing a nonzero gradient at $\mu^*$ along $\mathrm{span}(\ones) = \gA^H_{\text{bias}}$.

\item[-] \textbf{Angular drift under chiral $C_4$ symmetry (Modified Example~5)}: Suppose we weaken the symmetry of the standardized target $\tilde p$ in Example~5 from the full dihedral group $D_4$ to the pure rotational subgroup $H' = C_4$ (eliminating reflections) by setting $\tilde p(x) = q_0(R_\alpha^\top x)$ for an angle $\alpha \notin (\pi/4)\Z$. Because $C_4 \subset \SO(2)$ is abelian, conjugation by $r_k \in C_4$ acts trivially on $\SO(2)$, so $\Ad_{r_k} = \mathrm{id}$ and the trivial subspace expands from $\{0\}$ to the entire Lie algebra: $\gA^{C_4}_{\text{bias}} = \gA = \mathfrak{so}(2) \cong \R$. With reflections removed, symmetry imposes no constraint on the rotational gradient, and the variational minimizer shifts to $\theta^* + \alpha \pmod{\pi/2}$.
\end{itemize}

\paragraph{Infinitesimal generators and the Stein moment-matching condition.}
To characterize precisely when $\nabla_G \cL(g^*) = 0$ in the presence of a non-trivial bias subspace $\gA^H_{\text{bias}}$, we express the directional derivative of the pullback loss $J(u) = \KL(q_{\exp(u)} \| \tilde p)$ at $u = 0$ in terms of the log-density potentials of $q_0$ and $\tilde p$. Write
\[
q_0(z) \propto \exp(-\psi(z))\, , \qquad \tilde p(z) \propto \exp(-\varphi(z))\, ,
\]
where $\psi, \varphi : \cX \to \R$ are smooth $H$-invariant potential functions. Describing how an infinitesimal transport $u \in \gA$ deforms points $z \in \cX$ requires the fundamental vector field (or \emph{velocity field}) induced by the group action:

\begin{definition}[Infinitesimal generator / velocity field] \label{def:velocity-field}
For each tangent direction $u \in \gA$, the \emph{velocity field} (or \emph{infinitesimal generator}) $V_u: \cX \to T\cX$ is the vector field on $\cX$ assigning to each point $x \in \cX$ its instantaneous velocity under the flow $t \mapsto \exp(tu) \cdot x$:
\[
V_u(x) \coloneqq \left. \frac{\ud}{\ud t} \right|_{t=0} (\exp(tu) \cdot x)\, .
\]
\end{definition}

For our primary transport groups, the velocity field $V_u(x)$ takes an immediate closed form:
\begin{itemize}[leftmargin=*]
\item[-] \textbf{Additive translations} ($G = (\R^d, +)$, $\exp(tu) \cdot x = x + tu$): $V_u(x) = u$ (a constant vector field on $\R^d$).
\item[-] \textbf{Linear matrix transformations} ($G = \GL(d)$, $\exp(tU) \cdot x = e^{tU}x$): $V_U(x) = Ux$ (a linear vector field on $\R^d$).
\item[-] \textbf{Rotations on the sphere} ($G = \SO(3)$, $\gA \cong \R^3$): $V_\omega(x) = \omega \times x$ (the tangential rotational velocity field on $\cS^2$).
\end{itemize}

When a random variable $Z \sim q_0$ is pushed forward along the flow $x(t) = \exp(tu) \cdot Z$, the KL divergence $J(tu) = \KL(q_{\exp(tu)} \| \tilde p) = \E_{Z \sim q_0}[\log q_0(Z) - \log \tilde p(\exp(tu)\cdot Z) - \log|\det J_{\exp(tu)}(Z)|]$ varies at first order through two distinct mechanisms:
\begin{enumerate}[label=(\roman*), leftmargin=*]
  \item \emph{Cross-entropy / target potential variation}: Differentiating $-\log \tilde p(\exp(tu) \cdot Z) = \varphi(\exp(tu) \cdot Z) + \text{const}$ at $t = 0$ yields the directional derivative of the target potential along the velocity field, $\E_{q_0}[\langle \nabla \varphi(Z), V_u(Z) \rangle]$.
  \item \emph{Entropic / Jacobian volume variation}: By Jacobi's formula, differentiating the log-Jacobian determinant $-\log|\det J_{\exp(tu)}(Z)|$ at $t = 0$ yields the negative divergence of the velocity field, $-\E_{q_0}[\mathrm{div}(V_u)(Z)]$. Applying integration by parts (Stein's identity \citep{stein1981estimation}; cf.\ \citet{hyvarinen2005estimation}) against the density $q_0(z) \propto \exp(-\psi(z))$ re-expresses this volume expansion rate directly in terms of the base potential score:
  \[
  \E_{Z \sim q_0}\big[\mathrm{div}(V_u)(Z)\big] = -\E_{Z \sim q_0}\big[\langle \nabla \log q_0(Z), V_u(Z) \rangle\big] = \E_{Z \sim q_0}\big[\langle \nabla \psi(Z), V_u(Z) \rangle\big]\, .
  \]
\end{enumerate}
Combining these two terms yields our second main theorem, proven in detail in Appendix~\ref{app:moment-matching}.

\begin{theorem}[Stein moment-matching condition for stationarity] \label{thm:moment-matching}
Under the setting of Definition~\ref{def:vi-symmetry}, assume that $q_0 \propto \exp(-\psi)$ and $\tilde p \propto \exp(-\varphi)$ have differentiable (or locally Lipschitz and weakly differentiable a.e.) log-potentials $\psi$ and $\varphi$ satisfying the integrability condition (locally uniformly for small $|t|$)
$\E_{z \sim q_0}\!\big[\|\nabla \psi(z)\|\|V_u(z)\| + \|\nabla \varphi(z)\|\|V_u(z)\| + |\mathrm{div}(V_u)(z)|\big] < \infty$ and the boundary condition
$\lim_{R\to\infty} \int_{\partial B_R} q_0(z)\|V_u(z)\|\,\ud S(z) = 0$
(which holds automatically when $\cX$ is compact, and on $\R^d$ whenever $q_0$ has sub-Gaussian or faster tail decay and $\|V_u(z)\|$ grows at most polynomially). Then for any $u \in \gA$, the directional derivative of the pullback loss $J(u) = \KL(q_{\exp(u)}\|\tilde p)$ at $u = 0$ is given by
\[
\langle \nabla_G \cL(g^*), u \rangle_\gA = \langle \nabla J(0), u \rangle_\gA = \E_{z \sim q_0} \left[ \langle \nabla \varphi(z), V_u(z) \rangle \right] - \E_{z \sim q_0} \left[ \langle \nabla \psi(z), V_u(z) \rangle \right]\, .
\]
In particular, since $\nabla_G \cL(g^*) \in \gA^H_{\emph{bias}}$ by Proposition~\ref{prop:decomposition}, the true transport parameter $g^*$ is a first-order critical point of $\cL$ ($\nabla_G \cL(g^*) = 0$) if and only if the following moment-matching condition holds for every basis direction $u \in \gA^H_{\emph{bias}}$:
\begin{equation}\label{eq:moment-matching}
\E_{z \sim q_0} \left[ \langle \nabla \psi(z), V_u(z) \rangle \right] = \E_{z \sim q_0} \left[ \langle \nabla \varphi(z), V_u(z) \rangle \right]\, .
\end{equation}
\end{theorem}

Statistically, the quantity $\langle \nabla \psi(z), V_u(z) \rangle = -\langle \nabla \log q_0(z), V_u(z) \rangle$ is the \emph{directional score function} of the base model along the transformation generator $V_u$ (for instance, the radial virial score $z^\top \nabla \psi(z)$ measuring overall dispersion when $V_u(z) = z$ is an isotropic scaling, or the summed score $\ones^\top \nabla \psi(z)$ measuring common location when $V_u(z) = \ones$ is a diagonal shift). Equation~\eqref{eq:moment-matching} therefore has a transparent statistical meaning: the true parameter $g^*$ is a stationary point of the variational objective if and only if the standardized target $\tilde p$ and the base distribution $q_0$ have matching expected directional scores under $q_0$ along every symmetry-preserving bias direction $u \in \gA^H_{\text{bias}}$.

\textbf{Remark (Constructive debiasing by base distribution calibration).} When the moment-matching identity~\eqref{eq:moment-matching} is violated, $\nabla_G \cL(g^*)$ has a nonzero component in $\gA^H_{\text{bias}}$, reflecting a calibration mismatch between the variational base $q_0$ and the standardized target $\tilde p$ along the symmetry-preserving subgroup $\exp(\gA^H_{\text{bias}})$. Crucially, Theorem~\ref{thm:moment-matching} provides a direct constructive procedure to eliminate this bias without altering the symmetry structure of the variational family: one recalibrates the base distribution by pushing $q_0$ forward along an element $g_\Omega = \exp(u_\Omega)$ with $u_\Omega \in \gA^H_{\text{bias}}$, defining the calibrated base
\[
q_\Omega \coloneqq (g_\Omega)_\sharp q_0\, .
\]
Because every $u_\Omega \in \gA^H_{\text{bias}} \subseteq \gA^H$ generates transformations in the centralizer $G^H$ commuting with $H$, the recalibrated base $q_\Omega$ remains strictly $H$-symmetric. Specifically, replacing $q_0 \propto \exp(-\psi)$ by $q_\Omega \propto \exp(-\psi_\Omega)$ in~\eqref{eq:moment-matching} and pulling back to $z' \sim q_0$ via $z = g_\Omega \cdot z'$ (using Stein's identity on the left-hand side) yields the calibration equation
\[
\E_{z' \sim q_0}\big[\mathrm{div}(V_u)(g_\Omega \cdot z')\big] = \E_{z' \sim q_0}\big[\langle \nabla \varphi(g_\Omega \cdot z'), V_u(g_\Omega \cdot z') \rangle\big]\, , \qquad \forall u \in \gA^H_{\text{bias}}\, ,
\]
which admits a solution $u_\Omega \in \gA^H_{\text{bias}}$ whenever the left-hand side lies in the range of the right-hand side along $\exp(\gA^H_{\text{bias}})$. Because $\dim(\gA^H_{\text{bias}})$ is typically very small (equal to $1$ in both the Covariance and Permutation examples), solving this condition reduces to a simple scalar algebraic equation for $u_\Omega$. When the functional form of the standardized potential $\varphi$ is known (for example, when approximating a known heavy-tailed elliptical family by a Gaussian family), $u_\Omega$ can be computed offline prior to optimization; when $\varphi$ is unknown, one can alternatively eliminate $\gA^H_{\text{bias}}$ by restricting the transport group $G$ to a subgroup $G_0 \subset G$ (structural gauge fixing), as discussed in Section~\ref{sec:discussion}.

\paragraph{Guiding examples: Explicit moment-matching equations and base calibration.}
We conclude this section by working out the moment-matching condition~\eqref{eq:moment-matching} explicitly in our three examples with non-vanishing directional derivatives (full derivations are provided in Appendix~\ref{app:examples}):

\begin{itemize}[leftmargin=*]
\item[-] \textbf{Covariance (Example~4)}: Here $\gA^H_{\text{bias}} = \mathrm{span}(I_d)$ is one-dimensional, and the single basis generator $u = I_d$ induces the radial scaling velocity field $V_{I_d}(z) = I_d z = z$. For a recalibrated isotropic Gaussian base $q_\Omega = \cN(0, \lambda^2 I_d)$ with potential $\psi_\lambda(z) = \|z\|_2^2 / (2\lambda^2)$, the left-hand side of~\eqref{eq:moment-matching} is the virial expectation
\[
\E_{z \sim q_\Omega}\big[\langle \nabla \psi_\lambda(z), z \rangle\big] = \frac{1}{\lambda^2} \E_{z \sim q_\Omega}\big[\|z\|_2^2\big] = d\, .
\]
When the standardized target is $\tilde p = \cN(0, \sigma^2 I_d)$ with $\varphi(z) = \|z\|_2^2 / (2\sigma^2)$, the right-hand side is $\sigma^{-2} \E_{q_\Omega}[\|z\|_2^2] = d \lambda^2 / \sigma^2$, so~\eqref{eq:moment-matching} holds if and only if $\lambda = \sigma$. More generally, for any non-Gaussian spherically symmetric target $\tilde p(z) \propto \exp(-\varphi(z))$ with radial potential $\varphi(z) = \rho(\|z\|_2)$, we have $\langle \nabla \varphi(z), z \rangle = \rho'(\|z\|_2)\|z\|_2$, and the exact scale calibration condition for $\lambda > 0$ reads
\[
\frac{1}{\lambda^2}\E_{z \sim \cN(0, \lambda^2 I_d)}\big[\|z\|_2^2\big] = d = \E_{z \sim \cN(0, \lambda^2 I_d)}\big[\rho'(\|z\|_2)\|z\|_2\big]\, .
\]
Choosing $\lambda$ to satisfy this scalar equation eliminates the radial scale bias at $A^*$, making $\hat A = A^*$ stationary and recovering the target scatter parameter $\Sigma^* = A^*(A^*)^\top$ (while the fitted Gaussian density itself has covariance $\lambda^2 \Sigma^*$).

\item[-] \textbf{Permutation (Example~6)}: Here $\gA^H_{\text{bias}} = \mathrm{span}(\ones)$ is one-dimensional, and the generator $u = \ones$ induces the constant velocity field $V_{\ones}(z) = \ones$. The moment-matching condition~\eqref{eq:moment-matching} therefore requires the expected sum of partial derivatives of the potentials to match under $q_0$:
\[
\E_{z \sim q_0}\big[\langle \nabla \psi(z), \ones \rangle\big] = \E_{z \sim q_0}\big[\langle \nabla \varphi(z), \ones \rangle\big]\, .
\]
For the product Laplace base $q_0(z) \propto \exp(-\kappa\|z\|_1)$, we have $\langle \nabla \psi(z), \ones \rangle = \kappa \sum_{i=1}^n \mathrm{sign}(z_i)$, whose expectation under $q_0$ is $0$. If $\tilde p(z) \propto \exp(-\kappa\|z - \beta\ones\|_1)$, shifting the base distribution along $\gA^H_{\text{bias}}$ by $u_\Omega = \lambda \ones$ to $q_\Omega(z) \propto \exp(-\kappa\|z - \lambda\ones\|_1)$ achieves $\E_{q_\Omega}[\langle \nabla \varphi(z), \ones \rangle] = 0$ precisely when $\lambda = \beta$.

\item[-] \textbf{Square vs.\ Chiral $C_4$ rotation (Example~5)}: For the quartic potential $\psi(z) = \kappa(x^4 + y^4)$ on $\R^2$ (written in polar coordinates $z = (r\cos\phi, r\sin\phi)$ as $\psi(r,\phi) = \kappa r^4(\tfrac{3}{4} + \tfrac{1}{4}\cos 4\phi)$) and the tilted target pullback $\varphi(z) = \psi(R_{-\alpha} z)$, the rotational generator $\dot{\theta} \in \mathfrak{so}(2) \cong \R$ induces the tangential velocity field $V_{\dot{\theta}}(x,y) = \dot{\theta}(-y, x)^\top = \dot{\theta} \frac{\partial}{\partial \phi}$. The directional derivatives of $\psi$ and $\varphi$ along $V_{\dot{\theta}}$ are
\[
\langle \nabla \psi(z), V_{\dot{\theta}}(z) \rangle = -\kappa \dot{\theta} r^4 \sin(4\phi)\, , \qquad \langle \nabla \varphi(z), V_{\dot{\theta}}(z) \rangle = -\kappa \dot{\theta} r^4 \sin\big(4(\phi - \alpha)\big)\, .
\]
Because $q_0$ is invariant under the reflection $\phi \mapsto -\phi$, the odd trigonometric moment vanishes ($\E_{q_0}[r^4 \sin(4\phi)] = 0$), and expanding $\sin(4(\phi - \alpha)) = \sin(4\phi)\cos(4\alpha) - \cos(4\phi)\sin(4\alpha)$ yields
\[
\langle \nabla J(0), \dot{\theta} \rangle_\gA = -\kappa M_4 \sin(4\alpha) \, \dot{\theta}\, , \qquad \text{where } M_4 \coloneqq -\E_{z \sim q_0}\big[r^4 \cos(4\phi)\big] > 0\, .
\]
Consequently, $\nabla J(0) = 0$ if and only if $\sin(4\alpha) = 0$, i.e., $\alpha \equiv 0 \pmod{\pi/4}$, which is precisely the condition that $\tilde p$ possesses the reflection symmetries of $D_4$ aligned with the axes or diagonals of $q_0$.
\end{itemize}

%% ===========================================================================
\section{Discussion and future directions} \label{sec:discussion}
%% ===========================================================================

We have developed a general Lie-theoretic framework characterizing how distributional symmetries of the target and base measures interact with transport-based variational families. By studying the conjugation action of the symmetry group $H$ on the transport group $G$ and its linearization---the adjoint representation $\Ad_H$ on the transport Lie algebra $\gA$---Theorem~\ref{thm:gradient-constraint} and Proposition~\ref{prop:decomposition} establish that the Riemannian gradient $\nabla_G \cL(g^*)$ decomposes across three mutually orthogonal subspaces: the identifiability kernel $\mathfrak{n}$ (where the objective is locally constant), the symmetry-breaking subspace $\gA^\perp$ (where the gradient is forced to vanish by symmetry), and the symmetry-preserving bias subspace $\gA^H_{\text{bias}} = \gA^H \cap \mathfrak{n}^\perp$ (where model misspecification can induce a nonzero gradient). Theorem~\ref{thm:moment-matching} complements this geometric decomposition by providing an exact Stein-type moment-matching condition characterizing when the gradient along $\gA^H_{\text{bias}}$ vanishes.

\paragraph{Addressing the conjecture of \citet{margossian2024variational}.}
Our results provide a precise algebraic answer to the conjecture of \citet{margossian2024variational} that each symmetry of a target distribution $p$ corresponds to a statistic that variational inference can recover under mild conditions. Within the class of Lie group transformation models, given a compact symmetry group $H$ of the standardized target $\tilde p$, any $H$-normalized transport group $G$ defines a smooth manifold of identifiable parameters $\cM = G / (G \cap H)$ whose tangent space at the true coset $[g^*]$ is canonically identified with $\mathfrak{n}^\perp = \gA^H_{\text{bias}} \oplus \gA^\perp$. The components of $[g^*]$ along the symmetry-breaking directions $\gA^\perp$ are automatically first-order stationary points of the variational objective $\KL(q_g \| p)$ for \emph{any} $H$-symmetric base distribution $q_0$, without requiring parametric specification ($p \in \cQ_G$). Whenever $\gA^H_{\text{bias}} = \{0\}$, the entire coset $[g^*] \in G/N$ is recovered at first order by symmetry alone.

\paragraph{Structural gauge fixing versus base calibration.}
When the bias subspace $\gA^H_{\text{bias}}$ is nontrivial, our framework highlights two distinct remedies for eliminating estimation bias:
\begin{enumerate}[label=(\roman*), leftmargin=*]
  \item \emph{Base calibration via moment matching (Theorem~\ref{thm:moment-matching})}: One retains the full transport group $G$ and calibrates the base distribution $q_0 \mapsto q_\Omega = (\exp(u_\Omega))_\sharp q_0$ along $u_\Omega \in \gA^H_{\text{bias}}$ so that the directional score moments~\eqref{eq:moment-matching} match, recovering the full parameter $g^*$ (modulo $N$).
  \item \emph{Structural gauge fixing (subgroup restriction)}: When calibrating $q_0$ is impractical because the radial profile of $\tilde p$ is unknown, one can instead restrict the variational family to a smaller, $H$-normalized Lie subgroup $G_0 \subset G$ whose Lie algebra $\gA_0 = T_e G_0$ is a Lie subalgebra ($[\gA_0, \gA_0] \subseteq \gA_0$) transverse to $\gA^H_{\text{bias}}$, i.e., $\gA_0 \cap \gA^H_{\text{bias}} = \{0\}$. In the restricted family $\cQ_{G_0}$, the bias subspace vanishes identically, so the projection $g_0^*$ of $g^*$ onto $G_0 / (G_0 \cap H)$ is guaranteed to be a stationary point for any $H$-symmetric base $q_0$ (whenever $g^*$ factors as $g^* = g_0^* g_{\text{bias}}^* n$ with $g_0^* \in G_0$, $g_{\text{bias}}^* \in G^H$, and $n \in N$, so that pulling $p$ back by $g_0^*$ preserves $H$-invariance), at the cost of no longer estimating the parameter component along $\gA^H_{\text{bias}}$.
\end{enumerate}
Covariance estimation (Example~4, $G = \GL(d)$, $H = \Or(d)$) provides an instructive comparison between these strategies. Here the bias subspace is the one-dimensional trace direction $\gA^H_{\text{bias}} = \mathrm{span}(I_d)$. From a Lie-theoretic standpoint, the canonical $H$-normalized subgroup of $\GL(d)$ that excises $\mathrm{span}(I_d)$ is the special linear group $G_0 = \SL(d) = \{A \in \GL(d) : \det A = 1\}$, whose Lie algebra $\mathfrak{sl}(d) = \{M \in \mathfrak{gl}(d) : \mathrm{tr}(M) = 0\}$ is precisely $\mathfrak{n} \oplus \gA^\perp = \mathfrak{so}(d) \oplus \mathrm{Sym}_0(d)$. Restricting transport to $\SL(d)$ yields the quotient manifold $\cM_0 = \SL(d)/\SO(d)$ of unit-determinant shape matrices ($\det \Sigma = 1$), recovering the exact covariance shape while leaving only the single scalar generalized variance $\det \Sigma^*$ unestimated. By comparison, \citet{margossian2024variational} optimize over the full unconstrained group $G = \GL(d)$ (obtaining a stationary point $\hat\Sigma = c\Sigma^*$ proportional to $\Sigma^*$) and eliminate the scale factor $c > 0$ post-hoc by normalizing $\hat\Sigma$ into the \emph{correlation matrix} ($\Sigma_{ii} = 1$ for $i = 1, \dots, d$). Our decomposition clarifies both post-hoc normalization and prior subgroup restriction: from the unconstrained optimum $\hat\Sigma = c\Sigma^*$, one can in fact extract the full unit-determinant shape matrix $(\det \hat\Sigma)^{-1/d}\hat\Sigma \in \SL(d)/\SO(d)$ (preserving the $d-1$ ratios of marginal variances that are discarded by the correlation matrix), while if one wishes to restrict the optimization domain a priori to a subgroup $G_0 \subset \GL(d)$, restricting to $\SL(d)$ preserves full $\Or(d)$ conjugation invariance whereas enforcing unit diagonal entries breaks it. Similarly, in the Permutation model (Example~6, $G = (\R^n, +)$, $H = \Sn$), restricting $G$ to the $(n-1)$-dimensional zero-sum hyperplane $G_0 = \{ \mu \in \R^n : \sum_{i=1}^n \mu_i = 0 \} = \gA^\perp$ is normalized by $\Sn$ and guarantees unbiased recovery of all centered coordinate offsets $\mu^* - \bar{\mu}^* \ones$.

\paragraph{From first-order stationarity to second-order and global optimality.}
Theorem~\ref{thm:gradient-constraint} is a first-order result: it guarantees that $\nabla_G \cL(g^*) = 0$ whenever $\gA^H_{\text{bias}} = \{0\}$ (or when moment matching holds), which means $g^*$ is a critical point of $\cL$. Whether $g^*$ is a strict local minimizer (modulo the flat coset fiber $N$) or a global minimizer depends on higher-order curvature properties of the target and base distributions. Specifically, $g^*$ defines a strict local minimum on the quotient manifold $\cM = G/N$ if the Riemannian Hessian $\mathrm{Hess}_G \cL(g^*) = \nabla^2 J(0)$ is strictly positive definite on the horizontal subspace $\mathfrak{n}^\perp = \gA^H_{\text{bias}} \oplus \gA^\perp$. When the target $p$ is strongly multimodal---for instance, in Example~1 when $p = \tfrac{1}{2}\cN(\mu^* + \delta, 1) + \tfrac{1}{2}\cN(\mu^* - \delta, 1)$ has two widely separated modes ($\delta \gg 1$) while $q_0 = \cN(0, 1)$ is unimodal with unit variance---the symmetry center $\mu^*$ remains a critical point by Theorem~\ref{thm:gradient-constraint}, but bifurcates from a global minimum into a local maximum (or a saddle point in higher dimensions) straddled by two symmetry-broken global minima near $\mu^* \pm \delta$. For translation families ($G = \R^d$), log-concavity of $p$ ensures convexity of $\cL(\mu)$ and hence global uniqueness of the stationary point $\mu^*$; extending such global uniqueness guarantees to non-abelian Lie groups via geodesic convexity on Riemannian quotient manifolds $\cM = G/N$ \citep{boumal2023introduction} is an important direction for future research.

\paragraph{Further extensions.}
Several additional directions merit investigation. First, in empirical applications, target distributions often satisfy $H$-conjugate symmetry only approximately; combining our representation-theoretic decomposition with perturbation bounds on $\nabla J(0)$ would quantify how an $\epsilon$-asymmetry in $\tilde p$ induces an $O(\epsilon)$ drift of the variational stationary point along $\gA^\perp$. Second, while Theorem~\ref{thm:gradient-constraint} already holds for any divergence functional that is invariant under diffeomorphisms of $\cX$ (including all $f$-divergences and R\'enyi $\alpha$-divergences, since steps~($\dagger$), (*), and~(**) in its proof rely solely on pushforward invariance), extending the Stein moment-matching formula of Theorem~\ref{thm:moment-matching} to general $f$-divergences and Wasserstein distances will yield analogous base-calibration recipes beyond the KL divergence. Finally, our assumption that $H$ acts on $\cX$ by linear isometries can be relaxed to arbitrary isometric actions of compact Lie groups on Riemannian manifolds $(\cX, \mathrm{vol}_{\cX})$: replacing the Euclidean gradient and divergence in Theorem~\ref{thm:moment-matching} by the Riemannian gradient $\nabla_{\cX}$ and Laplace--Beltrami divergence $\mathrm{div}_{\cX}$ leaves all our main theorems intact.

\bibliography{references}

@article{hall1965relationship,
  author = {Hall, W. J. and Wijsman, R. A. and Ghosh, J. K.},
  title = {The relationship between sufficiency and invariance with applications in sequential analysis},
  volume = {36},
  journal = {The Annals of Mathematical Statistics},
  number = {2},
  publisher = {Institute of Mathematical Statistics},
  pages = {575--614},
  year = {1965}
}

@article{margossian2025variational,
  title={Variational inference for uncertainty quantification: An analysis of trade-offs},
  author={Margossian, Charles C. and Pillaud-Vivien, Loucas and Saul, Lawrence K.},
  journal={Journal of Machine Learning Research},
  volume={26},
  number={202},
  pages={1--41},
  year={2025}
}

@inproceedings{ghafourpour2026variational,
  title={Variational inference via radial transport},
  author={Ghafourpour, Luca and Chewi, Sinho and Figalli, Alessio and Pooladian, Aram-Alexandre},
  booktitle={International Conference on Artificial Intelligence and Statistics},
  year={2026}
}

@article{marks2026symmetry,
  title={Symmetry guarantees statistic recovery in variational inference},
  author={Marks, Daniel and Paccagnan, Dario and van der Wilk, Mark},
  journal={arXiv preprint arXiv:2604.18310},
  year={2026}
}

@inproceedings{margossian2024variational,
  title={Variational inference in location-scale families: Exact recovery of the mean and correlation matrix},
  author={Margossian, Charles C. and Saul, Lawrence K.},
  booktitle={International Conference on Artificial Intelligence and Statistics},
  volume={258},
  pages={3466--3474},
  year={2025},
  organization={PMLR}
}

@book{boumal2023introduction,
  title={An Introduction to Optimization on Smooth Manifolds},
  author={Boumal, Nicolas},
  year={2023},
  publisher={Cambridge University Press}
}

@article{blei2017variational,
  title={Variational inference: A review for statisticians},
  author={Blei, David M. and Kucukelbir, Alp and McAuliffe, Jon D.},
  journal={Journal of the American Statistical Association},
  volume={112},
  number={518},
  pages={859--877},
  year={2017}
}

@inproceedings{ranganath2014black,
  title={Black box variational inference},
  author={Ranganath, Rajesh and Gerrish, Sean and Blei, David M.},
  booktitle={International Conference on Artificial Intelligence and Statistics},
  volume={33},
  pages={814--822},
  year={2014},
  organization={PMLR}
}

@book{absil2008optimization,
  title={Optimization Algorithms on Matrix Manifolds},
  author={Absil, Pierre-Antoine and Mahony, Robert and Sepulchre, Rodolphe},
  year={2008},
  publisher={Princeton University Press}
}

@book{hall2015lie,
  title={Lie Groups, {L}ie Algebras, and Representations: An Elementary Introduction},
  author={Hall, Brian C.},
  edition={2},
  year={2015},
  publisher={Springer},
  series={Graduate Texts in Mathematics},
  volume={222}
}

@book{fulton2013representation,
  title={Representation Theory: A First Course},
  author={Fulton, William and Harris, Joe},
  year={1991},
  publisher={Springer},
  series={Graduate Texts in Mathematics},
  volume={129}
}

@book{knapp2002lie,
  title={Lie Groups Beyond an Introduction},
  author={Knapp, Anthony W.},
  edition={2},
  year={2002},
  publisher={Birkh{\"a}user},
  series={Progress in Mathematics},
  volume={140}
}

@article{bondar1981amenability,
  title={Amenability: A survey for statistical applications of {H}unt--{S}tein and related conditions on groups},
  author={Bondar, J. V. and Milnes, P.},
  journal={Zeitschrift f{\"u}r Wahrscheinlichkeitstheorie und verwandte Gebiete},
  volume={57},
  number={1},
  pages={103--128},
  year={1981}
}

@inproceedings{kohler2020equivariant,
  title={Equivariant flows: exact likelihood generative learning for symmetric densities},
  author={K{\"o}hler, Jonas and Klein, Leon and No{\'e}, Frank},
  booktitle={International Conference on Machine Learning},
  pages={5361--5370},
  year={2020},
  organization={PMLR}
}

@article{margossian2025generalized,
  title={Generalized guarantees for variational inference in the presence of even and elliptical symmetry},
  author={Margossian, Charles C. and Rankin, Isaac E. and Saul, Lawrence K.},
  journal={arXiv preprint arXiv:2511.01064},
  year={2025}
}

@book{amari2016information,
  title={Information Geometry and Its Applications},
  author={Amari, Shun-ichi},
  year={2016},
  publisher={Springer},
  series={Applied Mathematical Sciences},
  volume={194}
}

@book{murray1993differential,
  title={Differential Geometry and Statistics},
  author={Murray, Michael K. and Rice, John W.},
  year={1993},
  publisher={Chapman \& Hall/CRC}
}

@inproceedings{cohen2016group,
  title={Group equivariant convolutional networks},
  author={Cohen, Taco and Welling, Max},
  booktitle={International Conference on Machine Learning},
  pages={2990--2999},
  year={2016}
}

@article{zellinger2026even,
  title={Even more guarantees for variational inference in the presence of symmetries},
  author={Zellinger, Lena and Vergari, Antonio},
  journal={arXiv preprint arXiv:2604.21407},
  year={2026}
}

@article{bronstein2021geometric,
  title={Geometric deep learning: Grids, groups, graphs, geodesics, and gauges},
  author={Bronstein, Michael M. and Bruna, Joan and Cohen, Taco and Veli{\v{c}}kovi{\'c}, Petar},
  journal={arXiv preprint arXiv:2104.13478},
  year={2021}
}

@article{weiler2019general,
  title={General {E}(2)-equivariant steerable {CNNs}},
  author={Weiler, Maurice and Cesa, Gabriele},
  journal={Advances in Neural Information Processing Systems},
  volume={32},
  year={2019}
}

@book{eaton1989group,
  title={Group Invariance Applications in Statistics},
  author={Eaton, Morris L.},
  year={1989},
  publisher={Institute of Mathematical Statistics and American Statistical Association}
}

@book{lehmann2005testing,
  title={Testing Statistical Hypotheses},
  author={Lehmann, Erich L. and Romano, Joseph P.},
  edition={3},
  year={2005},
  publisher={Springer}
}

@book{wijsman1990invariant,
  title={Invariant Measures on Groups and Their Use in Statistics},
  author={Wijsman, Robert A.},
  series={Institute of Mathematical Statistics Lecture Notes--Monograph Series},
  volume={14},
  publisher={Institute of Mathematical Statistics},
  year={1990}
}

@article{hyvarinen2005estimation,
  title={Estimation of non-normalized statistical models by score matching},
  author={Hyv{\"a}rinen, Aapo},
  journal={Journal of Machine Learning Research},
  volume={6},
  pages={695--709},
  year={2005}
}

@article{stein1981estimation,
  title={Estimation of the mean of a multivariate normal distribution},
  author={Stein, Charles M.},
  journal={The Annals of Statistics},
  volume={9},
  number={6},
  pages={1135--1151},
  year={1981}
}

@article{papamakarios2021normalizing,
  title={Normalizing flows for probabilistic modeling and inference},
  author={Papamakarios, George and Nalisnick, Eric and Rezende, Danilo Jimenez and Mohamed, Shakir and Lakshminarayanan, Balaji},
  journal={Journal of Machine Learning Research},
  volume={22},
  number={57},
  pages={1--64},
  year={2021}
}

@article{jordan1999introduction,
  title={An introduction to variational methods for graphical models},
  author={Jordan, Michael I. and Ghahramani, Zoubin and Jaakkola, Tommi S. and Saul, Lawrence K.},
  journal={Machine Learning},
  volume={37},
  number={2},
  pages={183--233},
  year={1999},
  publisher={Springer}
}

@article{wainwright2008graphical,
  title={Graphical models, exponential families, and variational inference},
  author={Wainwright, Martin J. and Jordan, Michael I.},
  journal={Foundations and Trends in Machine Learning},
  volume={1},
  number={1--2},
  pages={1--305},
  year={2008},
  publisher={Now Publishers}
}

@article{kobyzev2020normalizing,
  title={Normalizing flows: An introduction and review of current methods},
  author={Kobyzev, Ivan and Prince, Simon J. D. and Brubaker, Marcus A.},
  journal={IEEE Transactions on Pattern Analysis and Machine Intelligence},
  volume={43},
  number={11},
  pages={3964--3979},
  year={2021}
}

@inproceedings{rezende2015variational,
  title={Variational inference with normalizing flows},
  author={Rezende, Danilo Jimenez and Mohamed, Shakir},
  booktitle={International Conference on Machine Learning},
  pages={1530--1538},
  year={2015}
}

\newpage

\appendix

\section{Proof of Theorem~\ref{thm:moment-matching}} \label{app:moment-matching}

We provide the complete proof of the Stein moment-matching condition (Theorem~\ref{thm:moment-matching}). Recall that $q_0(z) = Z_0^{-1} \exp(-\psi(z))$ and $\tilde p(z) = Z_{\tilde p}^{-1} \exp(-\varphi(z))$ have differentiable potentials $\psi, \varphi : \cX \to \R$ satisfying the integrability and boundary conditions stated in Theorem~\ref{thm:moment-matching}.

\begin{proof}[Proof of Theorem~\ref{thm:moment-matching}]
Recall from the proof of Theorem~\ref{thm:gradient-constraint} that the pullback objective around $g^*$ is given by $J(u) = \cL(g^* \exp(u)) = \KL(q_{\exp(u)} \| \tilde p)$ for $u \in \gA$. We compute the directional derivative $\langle \nabla_G \cL(g^*), u \rangle_\gA = \langle \nabla J(0), u \rangle_\gA = \left.\frac{\ud}{\ud t}\right|_{t=0} J(tu)$ for an arbitrary tangent vector $u \in \gA$.

\textbf{Step 1: Pullback to the fixed base measure $q_0$.} Under the diffeomorphism $x = \exp(tu) \cdot z$ with $z \sim q_0$, the pushforward density $q_{\exp(tu)}$ satisfies $q_{\exp(tu)}(\exp(tu) \cdot z) \, |\det J_{\exp(tu)}(z)| = q_0(z)$. Substituting the change of variables $x = \exp(tu) \cdot z$ (with $\ud x = |\det J_{\exp(tu)}(z)| \, \ud z$) into the definition of the KL divergence gives
\[
J(tu) = \int_{\cX} q_0(z) \Big[ \log q_0(z) - \log \tilde p\big(\exp(tu) \cdot z\big) - \log \big|\det J_{\exp(tu)}(z)\big| \Big] \, \ud z\, .
\]
Note that the first integrand term $\log q_0(z)$ is independent of $t$ and therefore does not contribute to the derivative with respect to $t$.

\textbf{Step 2: Differentiation at $t = 0$.} Differentiating under the integral sign at $t = 0$ (justified by dominated convergence under a standard integrable local dominating bound on $t \in (-\epsilon, \epsilon)$ together with the integrability assumption of Theorem~\ref{thm:moment-matching}), we obtain
\begin{align}
\left.\frac{\ud}{\ud t}\right|_{t=0} J(tu) &= -\int_{\cX} q_0(z) \left[ \left.\frac{\ud}{\ud t}\right|_{t=0} \log \tilde p\big(\exp(tu) \cdot z\big) + \left.\frac{\ud}{\ud t}\right|_{t=0} \log \big|\det J_{\exp(tu)}(z)\big| \right] \ud z\, .
\end{align}
For the first term, since $\log \tilde p(x) = -\varphi(x) - \log Z_{\tilde p}$ and $\left.\frac{\ud}{\ud t}\right|_{t=0} (\exp(tu) \cdot z) = V_u(z)$ by Definition~\ref{def:velocity-field}, the chain rule yields
\[
\left.\frac{\ud}{\ud t}\right|_{t=0} \log \tilde p\big(\exp(tu) \cdot z\big) = \big\langle \nabla \log \tilde p(z), V_u(z) \big\rangle = -\big\langle \nabla \varphi(z), V_u(z) \big\rangle\, .
\]
For the second term, the flow map $\phi_t(z) = \exp(tu) \cdot z$ satisfies $\phi_0(z) = z$ (with $J_{\phi_0}(z) = I_d$) and $\left.\frac{\ud}{\ud t}\right|_{t=0} J_{\phi_t}(z) = \nabla V_u(z)$. Applying Jacobi's formula for the derivative of a log-determinant gives the divergence of the velocity field $V_u$:
\[
\left.\frac{\ud}{\ud t}\right|_{t=0} \log \big|\det J_{\exp(tu)}(z)\big| = \mathrm{tr}\big(\nabla V_u(z)\big) = \mathrm{div}(V_u)(z) = \sum_{i=1}^d \frac{\partial (V_u)_i}{\partial z_i}(z)\, .
\]

\textbf{Step 3: Integration by parts (Stein's identity).} Using the product rule $\mathrm{div}(q_0 V_u) = q_0 \, \mathrm{div}(V_u) + \langle \nabla q_0, V_u \rangle$ and the divergence theorem on $\cX$, the boundary integral vanishes by our boundary condition, and since $\nabla q_0(z) = q_0(z) \nabla \log q_0(z) = -q_0(z) \nabla \psi(z)$, we obtain
\begin{align}
\int_{\cX} q_0(z)\, \mathrm{div}(V_u)(z)\, \ud z &= -\int_{\cX} \big\langle \nabla q_0(z), V_u(z) \big\rangle\, \ud z \\
&= -\int_{\cX} q_0(z) \big\langle \nabla \log q_0(z), V_u(z) \big\rangle\, \ud z = \int_{\cX} q_0(z) \big\langle \nabla \psi(z), V_u(z) \big\rangle\, \ud z\, .
\end{align}

\textbf{Step 4: Combining both contributions.} Substituting the results of Steps~2 and~3 into the derivative of $J(tu)$ yields
\[
\langle \nabla_G \cL(g^*), u \rangle_\gA = \langle \nabla J(0), u \rangle_\gA = \E_{z \sim q_0}\Big[ \big\langle \nabla \varphi(z), V_u(z)\big\rangle \Big] - \E_{z \sim q_0}\Big[ \big\langle \nabla \psi(z), V_u(z)\big\rangle \Big]\, .
\]
Finally, by Proposition~\ref{prop:decomposition}, $\nabla_G \cL(g^*)$ already lies in $\gA^H_{\text{bias}}$ (its components along $\mathfrak{n}$ and $\gA^\perp$ vanish identically). Thus $\nabla_G \cL(g^*) = 0$ if and only if $\langle \nabla_G \cL(g^*), u \rangle_\gA = 0$ for all $u \in \gA^H_{\text{bias}}$, which is equivalent to~\eqref{eq:moment-matching}.
\end{proof}

\section{$\Ad_H$-invariance of $\mathfrak{n}$ and subspace decomposition details} \label{app:centrality}

The orthogonal decomposition $\gA = \mathfrak{n} \oplus \gA^H_{\text{bias}} \oplus \gA^\perp$ in Proposition~\ref{prop:decomposition} uses the fact that the identifiability kernel $\mathfrak{n} = T_e N$ is an $\Ad_H$-invariant subspace of $\gA$. We record the detailed verification here and display the explicit decomposition $\mathfrak{n} = (\mathfrak{n} \cap \gA^H) \oplus (\mathfrak{n} \cap (\gA^H)^\perp)$ across all six guiding examples.

\subsection{Proof of $\Ad_H$-invariance of $\mathfrak{n}$}

Recall that $N = G \cap H$ is the intersection of the transport group $G$ and the symmetry group $H$ inside the total group $K = GH$. For any $n \in N$ and $h \in H$, we have $\Phi_h(n) = h \circ n \circ h^{-1} \in G$ because $G$ is normalized by $H$ (Definition~\ref{def:vi-symmetry}), and simultaneously $h \circ n \circ h^{-1} \in H$ because $n \in H$ and $H$ is a group. Hence $\Phi_h(n) \in G \cap H = N$ for all $n \in N$ and $h \in H$.

Now let $u \in \mathfrak{n}$. By definition of $\mathfrak{n}$, the one-parameter subgroup $t \mapsto \exp(tu)$ lies entirely in $N$ for all $t \in \R$. Applying $\Phi_h$, the smooth curve $t \mapsto \Phi_h(\exp(tu))$ also lies entirely in $N$ for all $t \in \R$, and passes through $e_G$ at $t = 0$. Its velocity vector at $t = 0$ therefore belongs to the tangent space $T_e N = \mathfrak{n}$:
\[
\Ad_h \cdot u = \left.\frac{\ud}{\ud t}\right|_{t=0} \Phi_h(\exp(tu)) \in \mathfrak{n}\, .
\]
Thus $\Ad_h(\mathfrak{n}) = \mathfrak{n}$ for every $h \in H$. Because the inner product $\langle \cdot, \cdot \rangle_\gA$ is $\Ad_H$-invariant, the orthogonal averaging projection $P_{\gA^H} = \int_H \Ad_h \, \ud\mathrm{Haar}(h)$ onto the trivial representation $\gA^H$ maps $\mathfrak{n}$ into $\mathfrak{n} \cap \gA^H$, and the complementary projection $I - P_{\gA^H}$ maps $\mathfrak{n}$ into $\mathfrak{n} \cap (\gA^H)^\perp$. Consequently, $\mathfrak{n}$ splits orthogonally as
\[
\mathfrak{n} = (\mathfrak{n} \cap \gA^H) \;\oplus\; (\mathfrak{n} \cap (\gA^H)^\perp)\, ,
\]
which ensures that $\mathfrak{n}$, $\gA^H_{\text{bias}} = \gA^H \cap \mathfrak{n}^\perp$, and $\gA^\perp = (\gA^H)^\perp \cap \mathfrak{n}^\perp$ are mutually orthogonal and sum to $\gA$.

\subsection{Explicit subspace verification for the six guiding examples}

\begin{itemize}[leftmargin=*]
\item[-] \textbf{Separable regime (Examples~1, 2, and~6)}: Here $N = \{e\}$, so $\mathfrak{n} = \{0\}$ and both $\mathfrak{n} \cap \gA^H = \{0\}$ and $\mathfrak{n} \cap (\gA^H)^\perp = \{0\}$ are trivial. Consequently, $\gA^H_{\text{bias}} = \gA^H$ and $\gA^\perp = (\gA^H)^\perp$.
\item[-] \textbf{Sphere (Example~3)}: Here $N = H = \SO(2)$ is abelian, so the adjoint action of $\SO(2)$ on its own Lie algebra $\mathfrak{n} = \mathrm{span}(e_3)$ is trivial, giving $\mathfrak{n} = \gA^H = \mathrm{span}(e_3)$. Thus $\mathfrak{n} \cap \gA^H = \mathfrak{n}$ and $\mathfrak{n} \cap (\gA^H)^\perp = \{0\}$: the entire identifiability kernel lies inside the trivial subspace $\gA^H$, leaving $\gA^H_{\text{bias}} = \gA^H \cap \mathfrak{n}^\perp = \{0\}$.
\item[-] \textbf{Covariance (Example~4)}: Here $N = H = \Or(d)$ has Lie algebra $\mathfrak{n} = \mathfrak{so}(d)$ (skew-symmetric matrices), whereas $\gA^H = \mathrm{span}(I_d)$ (isotropic scalar matrices). Because the only matrix that is both symmetric (a multiple of $I_d$) and skew-symmetric is the zero matrix, we have $\mathfrak{n} \cap \gA^H = \{0\}$ and $\mathfrak{n} \cap (\gA^H)^\perp = \mathfrak{so}(d)$: for $d \ge 2$, the entire identifiability kernel $\mathfrak{so}(d)$ lies in the non-trivial representation $(\gA^H)^\perp$. Consequently, $\gA^H_{\text{bias}} = \gA^H = \mathrm{span}(I_d)$, and $\gA^\perp = (\gA^H)^\perp \cap \mathfrak{so}(d)^\perp = \mathrm{Sym}_0(d)$ is the space of traceless symmetric matrices.
\item[-] \textbf{Square (Example~5)}: Here $N = C_4$ is a finite group of order $4$, so its Lie algebra is zero-dimensional ($\mathfrak{n} = \{0\}$), and the decomposition of $\mathfrak{n}$ is again trivial.
\end{itemize}

\section{Detailed derivations for the guiding examples} \label{app:examples}

\subsection{Covariance example: Explicit scale bias and base calibration}

Consider Example~4 ($G = \GL(d)$, $H = \Or(d)$) with standard Gaussian base $q_0 = \cN(0, I_d)$ and standardized target pullback $\tilde p = \cN(0, \sigma^2 I_d)$ for some $\sigma > 0$. The Lie algebra is $\gA = \mathfrak{gl}(d)$, and by Proposition~\ref{prop:decomposition} the bias subspace is $\gA^H_{\text{bias}} = \mathrm{span}(I_d)$.

For the basis generator $u = I_d$, the velocity field is $V_{I_d}(z) = I_d z = z$. The log-density potentials (up to additive constants) are $\psi(z) = \tfrac{1}{2}\|z\|_2^2$ and $\varphi(z) = \tfrac{1}{2\sigma^2}\|z\|_2^2$, with gradients $\nabla \psi(z) = z$ and $\nabla \varphi(z) = \sigma^{-2} z$. Substituting these into Theorem~\ref{thm:moment-matching} gives the directional derivative of the pullback loss along $I_d$:
\begin{align}
\langle \nabla_G \cL(A^*), I_d \rangle_{\text{Frob}} &= \E_{z \sim q_0}\big[\langle \nabla \varphi(z), z\rangle\big] - \E_{z \sim q_0}\big[\langle \nabla \psi(z), z\rangle\big] \\
&= \frac{1}{\sigma^2} \E_{z \sim q_0}\big[\|z\|_2^2\big] - \E_{z \sim q_0}\big[\|z\|_2^2\big] = d\left(\frac{1}{\sigma^2} - 1\right)\, .
\end{align}
Because $\|I_d\|_{\text{Frob}}^2 = d$ and $\nabla_G \cL(A^*) \in \mathrm{span}(I_d)$, we obtain the exact Riemannian gradient at $A^*$:
\[
\nabla_G \cL(A^*) = \left(\frac{1}{\sigma^2} - 1\right) I_d\, .
\]
This vanishes if and only if $\sigma^2 = 1$. When $\sigma^2 \neq 1$, recalibrating the base distribution along $\exp(\gA^H_{\text{bias}}) = \{\lambda I_d : \lambda > 0\}$ by choosing $g_\Omega = \sigma I_d$ produces the calibrated base $q_\Omega = (g_\Omega)_\sharp q_0 = \cN(0, \sigma^2 I_d)$ with potential $\psi_\Omega(z) = \tfrac{1}{2\sigma^2}\|z\|_2^2$. The moment-matching condition~\eqref{eq:moment-matching} is then satisfied identically, and $\nabla_G \cL(A^*) = 0$.

\subsection{Square example: Polar coordinate derivation of the angular gradient}

Consider Example~5 ($G = \SO(2)$, $H = D_4$) with quartic base distribution $q_0(z) \propto \exp(-\kappa(x^4 + y^4))$ on $\R^2$, where $z = (x,y)^\top$ and $\kappa > 0$. Passing to polar coordinates $x = r\cos\phi$, $y = r\sin\phi$ and using the identity $\cos^4\phi + \sin^4\phi = (\cos^2\phi + \sin^2\phi)^2 - \tfrac{1}{2}\sin^2(2\phi) = \tfrac{3}{4} + \tfrac{1}{4}\cos(4\phi)$, the base potential takes the form
\[
\psi(r, \phi) = \kappa(r^4 \cos^4\phi + r^4 \sin^4\phi) = \kappa r^4 \left(\frac{3}{4} + \frac{1}{4}\cos(4\phi)\right)\, .
\]
For an angular velocity $u = \dot{\theta} \in \mathfrak{so}(2) \cong \R$, the flow $\exp(t\dot{\theta})\cdot z = R_{t\dot{\theta}} z$ has velocity field $V_{\dot{\theta}}(x,y) = \dot{\theta}(-y, x)^\top$. Consequently, the directional derivative of $\psi$ along $V_{\dot{\theta}}$ is
\begin{align}
T_0(z) &\coloneqq \big\langle \nabla \psi(z), V_{\dot{\theta}}(z)\big\rangle = 4\kappa \dot{\theta} \left(-x^3 y + x y^3\right) = \dot{\theta} \frac{\partial \psi}{\partial \phi}(r, \phi) = -\kappa \dot{\theta} r^4 \sin(4\phi)\, .
\end{align}
Now suppose the standardized target pullback is rotated by an angle $\alpha$, so that $\varphi(z) = \psi(R_{-\alpha} z) = \kappa r^4\big(\tfrac{3}{4} + \tfrac{1}{4}\cos(4(\phi - \alpha))\big)$. The directional derivative of $\varphi$ along $V_{\dot{\theta}}$ is similarly
\[
T_\alpha(z) \coloneqq \big\langle \nabla \varphi(z), V_{\dot{\theta}}(z)\big\rangle = -\kappa \dot{\theta} r^4 \sin\big(4(\phi - \alpha)\big)\, .
\]
By Theorem~\ref{thm:moment-matching}, the directional derivative of the pullback loss $J$ at $0$ is
\begin{align}
\langle \nabla J(0), \dot{\theta} \rangle_\gA &= \E_{z \sim q_0}\big[T_\alpha(z) - T_0(z)\big] = -\kappa \dot{\theta} \, \E_{z \sim q_0}\Big[r^4\big(\sin(4(\phi-\alpha)) - \sin(4\phi)\big)\Big]\, .
\end{align}
Because the density $q_0(r,\phi) \propto \exp(-\kappa r^4(\tfrac{3}{4} + \tfrac{1}{4}\cos(4\phi)))$ is an even function of $\phi$ while $\sin(4\phi)$ is odd, we have $\E_{q_0}[r^4 \sin(4\phi)] = 0$. Expanding $\sin(4(\phi-\alpha)) = \sin(4\phi)\cos(4\alpha) - \cos(4\phi)\sin(4\alpha)$ therefore simplifies the gradient to
\[
\langle \nabla J(0), \dot{\theta} \rangle_\gA = \kappa \dot{\theta} \sin(4\alpha) \, \E_{z \sim q_0}\big[r^4 \cos(4\phi)\big] = -\kappa M_4 \sin(4\alpha) \, \dot{\theta}\, ,
\]
where $M_4 \coloneqq -\E_{z \sim q_0}[r^4 \cos(4\phi)]$. Note that $M_4$ is strictly positive for every $\kappa > 0$: for each fixed radius $r > 0$, the conditional angular density $p(\phi \mid r) \propto \exp(-\tfrac{1}{4}\kappa r^4 \cos(4\phi))$ places strictly greater mass on angles where $\cos(4\phi) < 0$ (near the diagonals $\phi = \pm \pi/4, \pm 3\pi/4$) than on angles where $\cos(4\phi) > 0$ (near the coordinate axes), ensuring $\E_{q_0}[r^4 \cos(4\phi)] < 0$. Thus $\nabla J(0) = 0$ if and only if $\sin(4\alpha) = 0$, i.e., $\alpha \in (\pi/4)\Z$.

\end{document}